\documentclass[
               ]{amsart}

\usepackage{float}
\usepackage
{graphicx}
\usepackage{
epsfig, epstopdf,
comment, tikz-cd, hyperref}

\newcommand{\ad}{{\operatorname{ad}}}
\newcommand{\End}{{\operatorname{End}}}
\newcommand{\Hom}{{\operatorname{Hom}}}

\newcommand{\image}{{\operatorname{Im}}}
\newcommand{\Ad}{{\operatorname{Ad}}}
\newcommand{\BHD}{{\operatorname{BHD}}}

\newcommand{\SL}{\operatorname{\rm{SL}}}
\newcommand{\BSL}{\operatorname{\rm{BSL}}}
\newcommand{\bx}{\boldsymbol{x}}
\newcommand{\by}{\boldsymbol{y}}

\theoremstyle{plain}
\newtheorem{theorem}{Theorem}[section]

\newtheorem{lemma}[theorem]{Lemma}
\newtheorem{proposition}[theorem]{Proposition}

\theoremstyle{definition}
\newtheorem{definition}[theorem]{Definition}

\theoremstyle{remark}

\allowdisplaybreaks

\begin{document}

\title[Quantized representations]{Quantized 
representations of knot groups}

\author[J. Murakami and R. van der Veen]{Jun Murakami and Roland van der Veen}

\thanks{The first author was  supported by JSPS KAKENHI Grant Number 17K187288.
The second author was supported by the Netherlands Organisation for Scientific Research.} 

\address[murakami@waseda.jp]{Jun Murakami, 
Waseda University,
Tokyo, Japan\vspace{-2mm}}

\address[roland.mathematics@gmail.com]{Roland van der Veen, University of Groningen, Groningen, the Netherlands}



\begin{abstract}
We propose a new non-commutative generalization of the representation variety and the character variety of
a knot group. Our strategy is to reformulate the construction of the algebra of functions on the space of representations
in terms of Hopf algebra objects in a braided category (braided Hopf algebra). The construction
works under the assumption that the algebra is braided commutative. The resulting knot invariant is a module with a coadjoint
action. Taking the coinvariants yields a new quantum character variety that may be thought of as an alternative to the skein module.
We give concrete examples for a few of the simplest knots and links.
\end{abstract}

\maketitle


\section{Introduction}\label{sec:intro}

The discovery of the Jones polynomial brought us a new method to study knots and links, but its relation to the geometric properties of the knot complement was unclear at that moment. After Witten's interpretation in terms of $\mathrm{SU}(2)$ Chern-Simons theory,  
R. Kashaev \cite{K} observed a precise relation between quantum invariants and the hyperbolic volume of the knot complement. This was reinterpreted as a relation between the colored Jones invariant and the hyperbolic volume by H. Murakami and the first author in \cite{MM}.  
Moreover, it was observed by Q. Chen and T. Yang in \cite{CY} that such relation also holds for the Witten-Reshetikhin-Turaev invariant of closed 3-manifolds.   
These relations between quantum invariants and hyperbolic volumes are not rigorously proved yet in general and are known as the volume conjecture.
In some sense, the volume conjecture means that the colored Jones invariants represent a quantization of the hyperbolic volume.  
Viewing the hyperbolic structure as a particular flat $\SL(2, {\mathbb C})$ connection, the above was given an interpretation in terms of topological quantum field theory with gauge group $\SL(2, {\mathbb C})$, see \cite{DG11}.
\par 
Once we got a relation like the volume conjecture, it is natural to think about quantization of other geometric properties.  For example if the knot is hyperbolic, its complement will be isometric to a quotient of hyperbolic space ${\mathbb H}^3/\Gamma$. The discrete subgroup  $\Gamma$ is isomorphic to the fundamental group of the knot complement (the knot group). Using a suitable quantization of the Lie group $\SL(2, {\mathbb C})$ and its discrete subgroup $\Gamma$, then the geometric structure of the complement of $K$ may be quantized.  More generally we construct a quantum deformation of the space of representations of the knot group into some linear algebraic group. 
\par
Our construction is not 3-dimensional but (2+1)-dimensional, as the construction of quantum invariants 
including the Jones polynomial and the Kontsevich invariant \cite{OhtsukiBook}.  
For a knot $K$, 
these invariants are obtained from a braid $b$ whose closure is isotopic to $K$.  
The braid $b$ is interpreted as an isotopy of a punctured disk, where the punctured disk is 2-dimensional and the deformation parameter is 1-dimensional.  
For the Jones polynomial, the braid group action is given by the quantum $R$-matrix, which comes from the monodromy matrix of conformal field theory.  
For the Kontsevich invariant, the braid group action is given by the Kontsevich's iterated integral.  
In both cases, we first consider such action of $b$, and then taking the `quantum trace' of these actions to get an invariant of $K$.  
\par
The starting point for this paper is the space of $G$ representation of the knot group,  $G$ is a linear algebraic group whose coordinate ring has a natural cocomutative Hopf algebra structure.    
Presenting the knot as a closed braid and interpreting the braid group action in terms of Hopf algebra we get a description of the representation space that is suitable for generalization. 
Replacing the coordinate ring of $G$ by a braided Hopf algebra and redoing the exact same construction while taking care of the braiding allows us to quantize the space of representations.  
To construct certain `trace', we need evaluation and coevaluation maps, which we do not know how to construct for our case with  Hopf algebras and braided Hopf algebras because they might be infinite dimensional.  Instead of taking a trace, 
we just take the $b$ invariant part of the algebra corresponding to the thickened punctured disk, and then show that it is independent of the choice of $b$.  
\par
We start by briefly recalling the construction of the space of representations of the knot group $\Gamma_K$ into a group $G$ that we aim to generalize/quantize in this work. 
The space of representations is described by an ideal in a tensor power of the
coordinate algebra $\mathbb{C}[G]$. The coordinate algebra is generated by the matrix entries and any presentation of $\Gamma_K$ allows us to express the relations as polynomial equations in these matrix entries. 
\par
This construction works for any finitely presented group and any affine algebraic group and is independent of the chosen presentation, see \cite[Proposition 8.2]{BH}.
However, it is not clear how to generalize this ideal  in a non-commutative deformation (i.e. quantizing) because one would need some way to order the variables that no longer commute. 

For a knot  $K$, presented as the closure of a braid $b$,  the Wirtinger presentation tells us all relations are given by conjugation.
Viewing the relations as equations on the matrix elements of our representation defines an ideal $I_b$ as follows.
To prepare our generalization to the non-commutative world we construct the submodule $I_b$ using the commutative Hopf algebra structure of the coordinate ring
${\mathbb C}[G]$:
\begin{align*}
&\text{$\Delta: {\mathbb C}[G] \to {\mathbb C}[G]^{\otimes2}$ with $\Delta(f)(a_1 \otimes a_2) = f(a_1 a_2)$}
\quad \text{(comultiplication)}, 
\\
 &\text{$S : {\mathbb C}[G] \to {\mathbb C}[G]$ with $S(f)(a) = f(a^{-1})$}
 \qquad\qquad\quad \text{(antipode)},
 \\
&\text{$\varepsilon: {\mathbb C}[G] \to \mathbb{C}$ with $\varepsilon(f)=f(e)$ ($e$ : the  identity of $G$)}
\quad \text{(counit)}.  
 \end{align*}
If the braid $b$ is a product of the standard generators $\sigma_1,\dots, \sigma_n$, 
say $b = \sigma_{i_{k}}^{\varepsilon_{k}}\sigma_{i_{k-1}}^{\varepsilon_{k-1}}\cdots$
$\sigma_{i_2}^{\varepsilon_2}\sigma_{i_1}^{\varepsilon_1}$,   
the ideal $I_b$ is generated by  
\[
\psi_b - \mathrm{id}
:
{\mathbb C}[G]^{\otimes n} \to {\mathbb C}[G]^{\otimes n}  
\qquad
(i=1, 2, \cdots, n)
\]
where $\psi_b$ is given by 
\[
\psi_b(\alpha_1\otimes \alpha_2\otimes \cdots\otimes \alpha_n) 
=
\left(\psi_{\sigma_{i_k}^{\varepsilon_k}}\circ\psi_{\sigma_{i_{k-1}}^{\varepsilon_{k-1}}}
\circ\cdots \circ\psi_{\sigma_{i_{2}}^{\varepsilon_{2}}}\circ\psi_{\sigma_{i_1}^{\varepsilon_1}}\right)(\alpha_1\otimes \alpha_2\otimes \cdots\otimes \alpha_n),
\]
and $\psi_{\sigma_{i}^{\pm1}}$ is given by
\begin{equation}
\begin{aligned}
\psi_{\sigma_{i}}(\alpha_1\otimes \cdots&\otimes \alpha_i\otimes \alpha_{i+1}\otimes \cdots\otimes \alpha_n)
\\
&=
(\alpha_1\otimes \cdots\otimes \alpha_{i+1}^{(3,2)}\otimes \alpha_i S(\alpha_{i+1}^{(3,1)}) \alpha_{i+1}^{(3,3)}\otimes \cdots\otimes \alpha_n),
\\
\psi_{\sigma_{i}^{-1}}(\alpha_1\otimes \cdots&\otimes \alpha_i\otimes \alpha_{i+1}\otimes \cdots\otimes \alpha_n)
\\
&=
(\alpha_1\otimes \cdots\otimes \alpha_{i}^{(3,1)} S(\alpha_{i}^{(3,3)}) \alpha_{i+1}\otimes \alpha_i^{(3,2)}\otimes \cdots\otimes \alpha_n).  
\end{aligned}
\label{eq:Ri}
\end{equation}
Here we use Sweedler's notation, i.e. the tensor $\alpha^{(3,1)} \otimes \alpha^{(3,2)} \otimes \alpha^{(3,3)}$ means $(\Delta\otimes id)(\Delta(\alpha))) = $ $\sum \alpha^{(3,1)} \otimes \alpha^{(3,2)} \otimes \alpha^{(3,3)}$. 

\begin{figure}[htb!]
\[
\begin{matrix}
\begin{matrix}
\begin{matrix}
\alpha \ \otimes \  \beta & & \alpha \ \otimes \  \beta
\\
\includegraphics[width=10mm,pagebox=cropbox,clip]{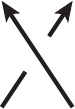}
&\quad&
\includegraphics[width=10mm,pagebox=cropbox,clip]{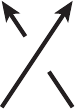} 
\\
\gamma \ \otimes \   \delta & &   \gamma \ \otimes \   \delta
\end{matrix}
&
\longrightarrow
\!\!\!\!\!
&
\begin{matrix}
\alpha \ \ \otimes \ \  \beta
& 
\alpha   \ \ \otimes \ \   \beta
\\
\includegraphics[keepaspectratio, scale=0.5]{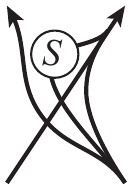}
&
\includegraphics[keepaspectratio, scale=0.5]{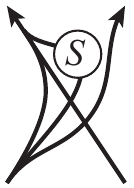} 
\\
\quad
\text{\footnotesize
$\begin{matrix}
\gamma\otimes\delta=\qquad\qquad \\
\beta^{(3,2)} \otimes \alpha\, S(\beta^{(3,1)})\beta^{(3,3)}
\end{matrix}$}
 &  
 \text{ 
 \footnotesize
 $\begin{matrix}
\gamma\otimes\delta=\qquad\qquad \\
\alpha^{(3,1)} S(\alpha^{(3,3)})\, \beta\otimes \alpha^{(3,2)}
\end{matrix}$}
\end{matrix}
\end{matrix}
\end{matrix}
\]
\caption{Dual  of the relations of the Wirtinger presentation at crossings and their diagrammatic presentations, read top to bottom.
}
\label{fig:dualWirtinger}
\end{figure}

\begin{figure}[htb!]
\begin{center}
$\begin{matrix}
\includegraphics[scale=0.9]{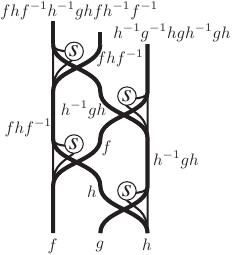}
\end{matrix}
\!\!\!\!\!\!\!\!\!\!\!\!\!\!\!\!\uparrow \pi_1
\qquad\ \ 
\begin{matrix}
\includegraphics[scale=0.9]{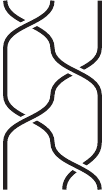}
\end{matrix}
\qquad
\mathbb{C}[G] \downarrow 
\!\!\!\!\!\!\!\!\!\!\!\!\!\!\!\!\!\!\!\!\!\!\!\!
\begin{matrix}
\includegraphics[scale=0.9]{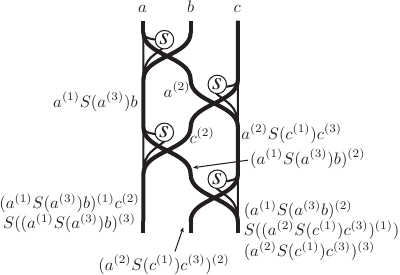}
\end{matrix}
$
\end{center}
\caption{The braid $(\sigma_2\sigma_1^{-1})^2$ interpreted as a Hopf diagram in the group algebra of $\pi_1$ (left, read bottom to top) or interpreted in $\mathbb{C}[G]$ (right, read top to bottom).}
\label{fig:Eightintro}
\end{figure}

As already mentioned each generator just acts by conjugation as in the Wirtinger presentation. 
A diagrammatic interpretation of \eqref{eq:Ri} is given in Figure \ref{fig:dualWirtinger}. The diagrams should be read top to bottom where each strand represents a copy of the algebra, the $Y$-shape represents the multiplication, the upside down $Y$ represents the coproduct and the $S$ represents the antipode, see also Figure \ref{fig:operation}.  In Figure \ref{fig:Eightintro} we showed what happens in the case of the braid $(\sigma_2\sigma_1^{-1})^2$ whose closure is the figure eight knot. Notice that reading the diagrams bottom to top and interpreting the $Y$-shape as the coproduct in the group algebra of the knot group recovers the corresponding Wirtinger presentation. 
\par
The construction of $I_b$ we sketched above works not only for ${\mathbb C}[\SL(2, {\mathbb C})]$ but also for any commutative Hopf algebra. 
Our main result is that it also works for braided commutative (braided) Hopf algebras. 
\par

\begin{figure}[htb!]
\[
\begin{matrix}
\begin{matrix}
\includegraphics[keepaspectratio, scale=0.7]{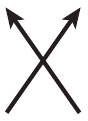}
\end{matrix}
&\ \ \longrightarrow\ \ &
\begin{matrix}
\includegraphics[scale=0.7]{cpositive.pdf}
\end{matrix}
\\
\text{\footnotesize flip} & & \text{\footnotesize braiding}
\end{matrix}
\]
\caption{Generalize the flip to the braiding}
\label{fig:braiding}
\end{figure}

A braided Hopf algebra $A$ is a generalization of a Hopf algebra where the braiding is used instead of the usual flip sending $x \otimes y$ to $y \otimes x$ as in Figure~\ref{fig:braiding}. Braided commutativity is a generalization of the commutativity property of usual Hopf algebras, which is given in Definition~\ref{def:braidedcomm}.     
\par
\begin{figure}[htb!]
\[
\begin{matrix}
\begin{matrix}
\includegraphics[scale=0.6]{cpositive.pdf}
\end{matrix}
\ \longrightarrow\ 
\begin{matrix}
\includegraphics[scale=0.6]{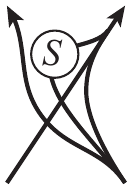}
\end{matrix}
\ \longrightarrow\ 
\begin{matrix}
\includegraphics[scale=0.6]{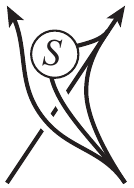}
\end{matrix},
&\qquad&
\begin{matrix}
\includegraphics[scale=0.6]{cnegative.pdf}
\end{matrix}
\ \longrightarrow\ 
\begin{matrix}
\includegraphics[scale=0.6]{cnegativediagram.pdf}
\end{matrix}
\ \longrightarrow\ 
\begin{matrix}
\includegraphics[scale=0.6]{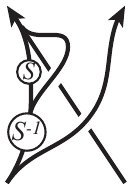}
\end{matrix}
\end{matrix}
\]
\caption{Relations of the braided Wirtinger presentation, read top to bottom.}
\label{fig:braidedWirtinger}
\end{figure}

To generalize the above construction of the ideal $I_b$ to get a space of $A$ representations, we modify the relation at the crossing as in Figure~\ref{fig:braidedWirtinger}.  
Our main result is to define a module $I_b$ and show that the quotient of $A^{\otimes n}$ divided by $I_b$ only depends on the knot $K$, see Theorem \ref{th:main}.  
  In the final example at the end of the paper we will return to the figure eight knot and show what our construction amounts to in this case.
\par
An important example of braided Hopf algebras is $\BSL(2)$, it is the braided  one-parameter deformation of the coordinate ring of  $\SL(2, {\mathbb C})$, see \cite{M91}.  
By applying the above construction, we get the space of $\BSL(2)$ representations which is a quantization of the $\SL(2, {\mathbb C})$ representation space of  $K$.  

Let $A_b{}^{A}$ be the $\Ad$ invariant subspace, i.e.
\[
A_b{}^{A}
=
\{x \in A^{\otimes n}/I_{d(b)} \mid {\rm Ad}(x) = x \otimes 1\}.
\]
We call $A_b{}^{A}$ the {\it quantum $A$ character variety} of $K$.  
If $A = \BSL(2)$, we also call it the {\it quantum $\SL(2)$ character variety}.  
Note that $A_b{}^{A}$ is not an algebra but an $\Ad$-comodule.
So the quantum character variety is not a variety in the usual sense.  

In the special case of $\SL(2)$, the quantum character variety we just defined seems to be equal to the skein module of the knot complement, which is often viewed as a quantization of the $\SL(2)$ character variety \cite{Le06}. 
\par
Our construction of quantum $\SL(2, {\mathbb C})$ character variety seems to be generated by quantum traces as in \cite{DL}.  
A more detailed discussion of our quantization of the quantum $\SL(2, {\mathbb C})$ character variety will appear in a forthcoming paper.  
More generally it seems plausible that our construction is related to the skein module defined for any ribbon category and any 3-manifold in \cite{GJS}, Definition 2.2.

\par
A similar definition of a quantum analogue of the character variety is also given by Habiro \cite{H}. It would also be interesting to
compare our quantization to the quantization based on ideal triangulations given in \cite{D13} and also with the quantization
procedure of \cite{BBJ18}. 

\par
This paper is organized as follows.  
In Section 2, we introduce the braided Hopf algebra with a focus on the braided commutative case.
We also introduce braided Hopf diagrams to explain morphisms between tensor powers of the braided Hopf algebra.
In Section 3, we construct representation of the braid group $B_n$ in $\End(A^{\otimes n})$ for any braided Hopf algebra $A$. 
Here we use the braided version of the Wirtinger presentation given in  Figure~\ref{fig:braidedWirtinger}.  

\begin{figure}[htb!]
\[
\begin{matrix}
A^{\otimes n}\otimes A^{\otimes n}
\\{}\\
\downarrow
\\{}\\
A^{\otimes n}
\end{matrix}
:
\qquad
\begin{matrix}
\begin{matrix}
\includegraphics[scale=0.7]{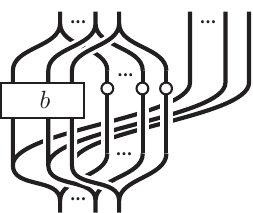}
\end{matrix}\ ,
&\qquad \qquad
&
\begin{matrix}
\includegraphics[scale=0.7]{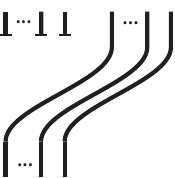}
\end{matrix}
\\
d_1\quad & & d_2\qquad
\end{matrix}
\]
\caption{Braided Hopf diagrams of $d_1$ and $d_2$, where  the space of $A$ representations $A_b = A^{\otimes n}/\image(d_1-d_2)$.}
\label{fig:db}
\end{figure}

In Section 4, we define the space of $A$ representations of a knot $K$ for any braided Hopf algebra $A$ satisfying braided commutativity.  
Let $b$ be a braid in $B_n$ whose closure $\widehat b$ is isotopic to $K$, and add $n$ strands to represent elements of the fundamental group twined to $b$ as in the Hopf diagram $d_1$ shown in Figure \ref{fig:db}.  
Let $I_b$ be the image of the map corresponding to $d_1-d_2$,  the space of $A$ representations is defined as  $A^{\otimes n}/I_b$.  
We show that this space only depends on the isotopy type of $\widehat b$.  

In Section 5, we apply the above construction to the trefoil knot, the Hopf link and the figure eight knot.  
%
\subsection*{Acknowlegement}
The authors would like to thank Professor Takefumi Nosaka for valuable comments. 
\section{Braided Hopf algebra and braided commutativity}
\subsection{Braided Hopf algebra}
\begin{figure}[htb]
\[
\begin{matrix}
\begin{matrix}
\includegraphics[scale=0.45]{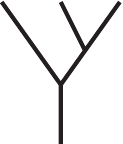}
\end{matrix}
=
\begin{matrix}
\includegraphics[scale=0.45]{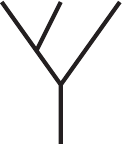}
\end{matrix}
\quad
\begin{matrix}
\includegraphics[scale=0.45]{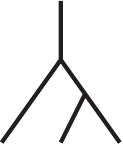}
\end{matrix}
=
\begin{matrix}
\includegraphics[scale=0.45]{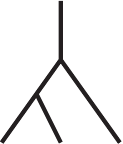}
\end{matrix}
\quad
\begin{matrix}
\includegraphics[scale=0.45]{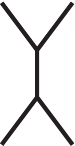}
\end{matrix}
=\ \ 
\begin{matrix}
\includegraphics[scale=0.45]{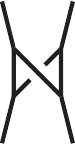}
\end{matrix}
\quad
\begin{matrix}
\includegraphics[scale=0.45]{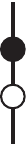}
\end{matrix}
\ \ =\ \ 
\begin{matrix}
\includegraphics[scale=0.45]{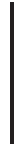}
\end{matrix}
\ \ =\ \ 
\begin{matrix}
\includegraphics[scale=0.45]{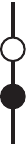}
\end{matrix}
\\
\begin{matrix}
\includegraphics[scale=0.45]{unitl}
\end{matrix}
\ =\ 
\begin{matrix}
\includegraphics[scale=0.45]{id}
\end{matrix}
\ =\ 
\begin{matrix}
\includegraphics[scale=0.45]{unitr}
\end{matrix}
\quad
\raisebox{2mm}{
$
\begin{matrix}
\includegraphics[scale=0.45,angle=180]{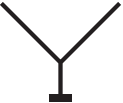}
\end{matrix}$
}
=
\raisebox{2mm}{
$\begin{matrix}
\includegraphics[scale=0.45,angle=180]{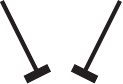}
\end{matrix}$
}
\quad
\begin{matrix}
\includegraphics[scale=0.45,angle=180]{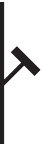}
\\
{}
\end{matrix}
\ =\ 
\begin{matrix}
\includegraphics[scale=0.45]{id}
\end{matrix}
\ =\ 
\raisebox{2mm}{
$\begin{matrix}
\includegraphics[scale=0.45,angle=180]{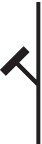}
\end{matrix}$
}
\quad
\begin{matrix}
\includegraphics[scale=0.45]{multunit}
\end{matrix}
=
\begin{matrix}
\includegraphics[scale=0.45]{unitmult}
\end{matrix}
\\
\ \ 
\begin{matrix}
\includegraphics[scale=0.45]{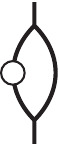}
\end{matrix}
\ \ =\ \ 
\begin{matrix}
\includegraphics[scale=0.45]{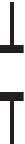}
\end{matrix}
\ \ =\ \ 
\begin{matrix}
\includegraphics[scale=0.45]{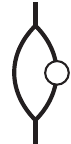}
\end{matrix}
\ \  \ \ 
\begin{matrix}
\includegraphics[scale=0.45]{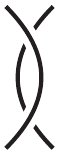}
\end{matrix}
\ \ =\ \ 
\begin{matrix}
\includegraphics[scale=0.45]{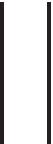}
\end{matrix}
\ \ =\ \ 
\begin{matrix}
\includegraphics[scale=0.45]{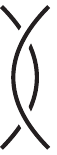}
\end{matrix}
\quad\quad
\begin{matrix}
\includegraphics[scale=0.45]{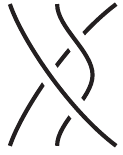}
\end{matrix}
\ =\ 
\begin{matrix}
\includegraphics[scale=0.45]{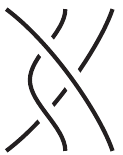}
\end{matrix}
\\
\ \ 
\begin{matrix}
\includegraphics[scale=0.45]{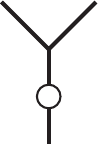}
\end{matrix}
=\ \ 
\begin{matrix}
\includegraphics[scale=0.45]{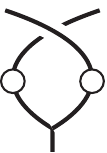}
\end{matrix}
\quad\quad
\raisebox{2mm}{
$\begin{matrix}
\includegraphics[scale=0.45,angle=180]{Sm1}
\end{matrix}$
}
=\ \ 
\raisebox{2mm}{
$\begin{matrix}
\includegraphics[scale=0.45,angle=180]{Sm2}
\end{matrix}$
}
\quad\quad
\begin{matrix}
\includegraphics[scale=0.45]{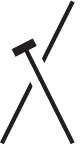}
\end{matrix}
=
\begin{matrix}
\includegraphics[scale=0.45]{unitpsi2}
\end{matrix}
\qquad
\raisebox{2mm}{
$\begin{matrix}
\includegraphics[scale=0.45,angle=180]{unitpsi1}
\end{matrix}$}
=
\raisebox{2mm}{
$\begin{matrix}
\includegraphics[scale=0.45,angle=180]{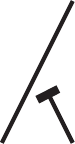}
\end{matrix}$
}
\\[24pt]
\quad
\includegraphics[scale=0.45]{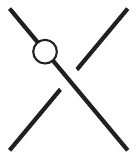}
\raisebox{0.5cm}{=}
\includegraphics[scale=0.45]{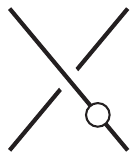}
\quad\quad
\includegraphics[scale=0.45]{multpsi1}
\raisebox{0.5cm}{=}
\includegraphics[scale=0.45]{multpsi2}
\quad\quad
\raisebox{1.2cm}{\includegraphics[scale=0.45,angle=180]{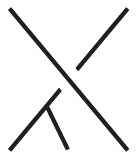}
\raisebox{-0.9cm}{=}
\includegraphics[scale=0.45,angle=180]{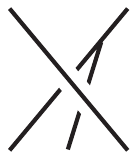}}
\end{matrix}
\]
\caption{The relations of a braided Hopf algebra, read from top to bottom.}
\label{fig:diagramrelations}
\end{figure}

A braided Hopf algebra is a version of a Hopf algebra having an extra operation called braiding.  
It may also be viewed as a Hopf object in a braided monoidal category. Such algebras are quite common in that
they can be produced from any quasi-triangular Hopf algebra by transmutation \cite{M02}. These structures also go by the name 
braided group.
\begin{definition}
An algebra $A$ over a field $k$ is called {\it a braided Hopf algebra} if it is equipped with following linear maps described by the diagrams in Figure~\ref{fig:operation}
satisfying the relations given in Figure~\ref{fig:diagramrelations}.   
\[
\begin{tabular}{rlrl}
$\text{ multiplication}$&
\!$\mu : A \otimes A \to A$,
&\!\!$\text{comultiplication}$&
$\!\Delta : A \to A \otimes A$,
\\
$\text{unit}$&\!\!
$\!1 : k \to A$,
&$\text{counit}$&
$\!\varepsilon : A \to k$,
\\
$\text{antipode}$&
$\!S : A \to A$,
&$\text{braiding}$&
$\!\Psi : A \otimes A \to A \otimes A$.
\end{tabular}
\]
\begin{figure}[htb]  
\[
\begin{matrix}
\begin{matrix}
\includegraphics[scale=0.45]{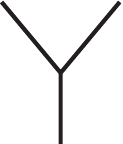}
\end{matrix}
&
\begin{matrix}
\raisebox{1.4cm}{\includegraphics[scale=0.45, angle=180]{mult.pdf}}
\end{matrix}
&
\begin{matrix}
\includegraphics[scale=0.45]{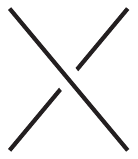}
\end{matrix}
&
\begin{matrix}
\includegraphics[scale=0.45]{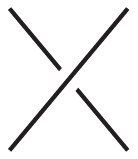}
\end{matrix}
\\
\text{\footnotesize multiplication} & \text{\footnotesize comultiplication $\Delta$} & \text{\footnotesize braiding $\Psi$} & \text{\footnotesize inverse braiding $\Psi^{-1}$}
\\{}\\
\begin{matrix}
\includegraphics[scale=0.45]{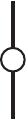}
\end{matrix}
&
\begin{matrix}
\includegraphics[scale=0.45]{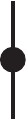}
\end{matrix}
&
\begin{matrix}
\includegraphics[scale=0.45]{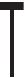}
\end{matrix}
&
\begin{matrix}\includegraphics[scale=0.45,angle=180]{unit.pdf}
\end{matrix}
\\
\text{\footnotesize antipode $S$} \quad&\qquad \text{\footnotesize $S^{-1}$} \qquad\quad&\ 
\text{\footnotesize unit $1$} \qquad&\quad \text{\footnotesize counit $\varepsilon$}\qquad
\end{matrix}
\]
\caption{The operations of the braided Hopf algebra $A$.}
\label{fig:operation}
\end{figure}
\label{def:bhopf}
\end{definition}
\begin{definition}
A diagram expressing a linear mapping from $A^{\otimes m}$ to $A^{\otimes n}$ built from a combination of the Hopf algebra operations given in Figure~\ref{fig:operation} is called 
{\it a braided Hopf diagram}.  
Let $\BHD(m, n)$ denote the set of braided Hopf diagrams expressing linear homomorphisms from $A^{\otimes m}$ to $A^{\otimes n}$.
\end{definition}


%
\subsection{Adjoint coaction}
A $k$-vector space $M$ is called a right $A$-comodule if 
there is a linear map
\[
\Delta : M \to M \otimes A
\]
satisfying the coassociativity
\[
(\Delta \otimes id)(\Delta)
=
(id \otimes \Delta)(\Delta).  
\]
Then $A$ itself is a right $A$-comodule with the following adjoint coaction $\ad : A \to A \otimes A$.  
\[
\ad(x)
=
(id \otimes \mu)(\Psi \otimes id)(S \otimes\Delta)\Delta(x),  
\]
\[
\begin{matrix}
\includegraphics[scale=0.6]{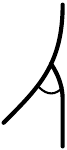}
\end{matrix}
\ \ =
\begin{matrix}
\includegraphics[scale=0.6]{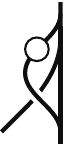}
\end{matrix}
\]
where $\mu : A \otimes A \to A$ is the multiplication of $A$, i.e. $\mu(x\otimes y) = x\, y$.  
\begin{proposition}[
C.f. \cite{M93b}, Proposition A.1]
\label{prop.Adrels}
Adjoint coaction satisfies the following relations.  
\begin{align}
&(id\otimes id \otimes \mu)(id \otimes \Psi \otimes id)
(\ad \otimes \ad)\Delta(x)
=
(\Delta \otimes id) \,\ad (x),  
\label{eq:adcomult}
\\
& (\ad \otimes id)\ad= (id \otimes \Delta)\ad.
\label{eq:addelta}
\\
& (\varepsilon \otimes id)\ad= 1\circ \varepsilon \qquad (id \otimes \varepsilon )\ad= id
\label{eq:adepsilon}
\end{align}
\[
\text{\eqref{eq:adcomult} : }\quad
\begin{matrix}
\includegraphics[scale=0.6]{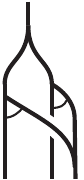}
\end{matrix}
\ =\ 
\begin{matrix}
\includegraphics[scale=0.6]{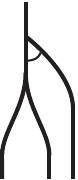}
\end{matrix}
\ 
,
\qquad
\text{\eqref{eq:addelta} : }\quad
\begin{matrix}
\includegraphics[scale=0.6]{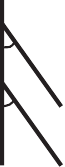}
\end{matrix}
\ =\ 
\begin{matrix}
\includegraphics[scale=0.6]{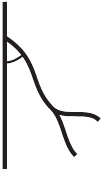}
\end{matrix}
\]
\label{prop:adcomm}
\end{proposition}
\begin{proof}
The relations \eqref{eq:adcomult} and \eqref{eq:addelta} are proved by the  graphical computation in Figure~\ref{fig:adcomm}.
The relations \eqref{eq:adepsilon} come from the properties of the unit $1$ and the counit $\varepsilon$.  
\begin{figure}[htb!]  
\[
\begin{matrix}
\begin{matrix}
\includegraphics[scale=0.7]{adcomult1}
\end{matrix}
\ \ =\ \ 
\begin{matrix}
\includegraphics[scale=0.7]{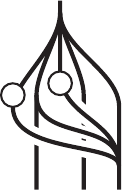}
\end{matrix}
\ \ =\ \ 
\begin{matrix}
\includegraphics[scale=0.7]{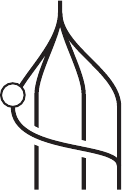}
\end{matrix}
\ \ =\ \ 
\begin{matrix}
\includegraphics[scale=0.7]{adcomult2}
\end{matrix}
\\
\text{\small\eqref{eq:adcomult} : Commutativity of adjoint and comultiplication.}
\\{}\\
\begin{matrix}
\includegraphics[scale=0.7]{addelta1}
\end{matrix}
\ =\ 
\begin{matrix}
\includegraphics[scale=0.7]{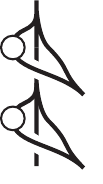}
\end{matrix}
\ =\ 
\begin{matrix}
\includegraphics[scale=0.7]{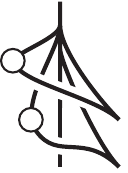}
\end{matrix}
\ =\ 
\begin{matrix}
\includegraphics[scale=0.7]{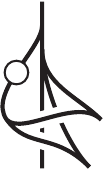}
\end{matrix}
\ =\ 
\begin{matrix}
\includegraphics[scale=0.7]{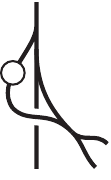}
\end{matrix}
\ =\ 
\begin{matrix}
\includegraphics[scale=0.7]{addelta6}
\end{matrix}
\\
\text{\small\eqref{eq:addelta} : Adjoint and multiplication.}
\end{matrix}
\]
\caption{Graphical proof of Proposition \ref{prop:adcomm}.}
\label{fig:adcomm}
\end{figure}
\end{proof}
\subsection{Braided commutativity}
\begin{figure}[htb!]
\[
\begin{matrix}
\includegraphics[scale=0.6]{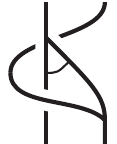}
\end{matrix}
\ =\ \ 
\begin{matrix}
\includegraphics[scale=0.6]{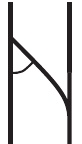}
\end{matrix}
\]
\caption{Graphical explanation for the braided commutativity.}
\label{fig:braidedcomm}
\end{figure}

\begin{figure}[htb!]
\[
\begin{matrix}
\includegraphics[scale=0.6]{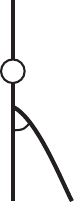}
\end{matrix}
\ =\ \ 
\begin{matrix}
\includegraphics[scale=0.6]{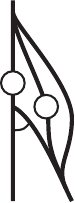}
\end{matrix}\ \ 
{\underset{\text{(bc)}}{=}}\ \ 
\begin{matrix}
\includegraphics[scale=0.6]{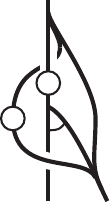}
\end{matrix}
\ \ =\ \ 
\begin{matrix}
\includegraphics[scale=0.6]{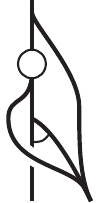}
\end{matrix}
\ \ =\ 
\begin{matrix}
\includegraphics[scale=0.6]{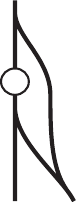}
\end{matrix}
\ =\ 
\begin{matrix}
\includegraphics[scale=0.6]{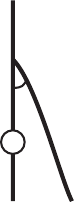}
\end{matrix}
\]
\caption{The adjoint $\ad$ commutes with $S$.}
\label{fig:adS}
\end{figure}

\begin{figure}[htb!]
\[
\begin{matrix}
\includegraphics[scale=0.7]{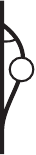}
\end{matrix}
\ =\ 
\begin{matrix}
\includegraphics[scale=0.7]{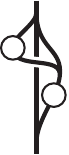}
\end{matrix}
\ =\ 
\begin{matrix}
\includegraphics[scale=0.7]{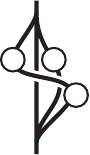}
\end{matrix}
\ =\ 
\begin{matrix}
\includegraphics[scale=0.7]{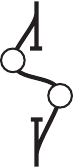}
\end{matrix}
\ =\ 
\begin{matrix}
\includegraphics[scale=0.7]{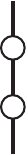}
\end{matrix}
\]
\caption{Proof for $\mu \circ (id \otimes S) \circ \ad = S^2$. }
\label{fig:adSS}
\end{figure}

We introduce the notion of the braided commutativity,  
which implies the compatibility of the adjoint coaction with respect to the multiplication $\mu$, the braiding $\Psi$, and the antipode $S$.  
\begin{definition}
The braided Hopf algebra $A$ is {\it braided commutative} if it satisfies
\begin{equation}
(id \otimes \mu)(\Psi\otimes id)(id\otimes \ad)  \Psi
=
(id \otimes \mu)(\ad \otimes id).
\label{eq:comutative}
\end{equation}
This relation is explained graphically in Figure \ref{fig:braidedcomm}.  
\label{def:braidedcomm}
\end{definition}
Braided commutativity was introduced in \cite{M93} and it is shown there that many interesting braided Hopf algebras
have this property. For example transmutation procedure always produces braided commutative braided Hopf algebras.
In the remainder of this section we assume that $A$ is braided commutative.
\begin{proposition}
The adjoint coaction commutes with the multiplication, i.e.
\begin{equation}
(\ad \circ \mu)
= 
(\mu \otimes \mu)(id \otimes \Psi \otimes id)(\ad \otimes \ad).
\label{eq:admult}
\end{equation}  
\[
\raisebox{-1cm}{\includegraphics[scale=0.6]{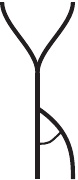}}
=\ \ 
\raisebox{-1cm}{\includegraphics[scale=0.6]{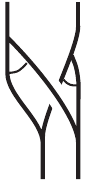}}
\]
\label{prop:admult}
\end{proposition}
\begin{proof}
The relation \eqref{eq:admult} is proved by the  graphical computation in Figure~\ref{fig:mult}.  
At the second to last equality, we use the braided comutativity.  
\begin{figure}[htb!]
\[
\begin{matrix}
\includegraphics[scale=0.6]{admult1}
\end{matrix}
=\ \ 
\begin{matrix}
\includegraphics[scale=0.6]{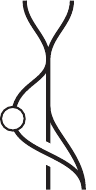}
\end{matrix}
=\ \ 
\begin{matrix}
\includegraphics[scale=0.6]{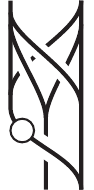}
\end{matrix}
\ \ =\ \ 
\begin{matrix}
\includegraphics[scale=0.6]{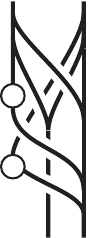}
\end{matrix}
\ \ =\ \ 
\begin{matrix}
\includegraphics[scale=0.6]{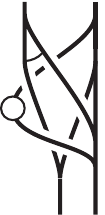}
\end{matrix}
{\underset{\text{(bc)}}{=}}
\begin{matrix}
\includegraphics[scale=0.6]{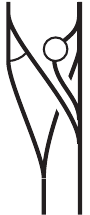}
\end{matrix}
\ \ =\ 
\begin{matrix}
\includegraphics[scale=0.6]{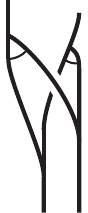}
\end{matrix}
\]
\caption{Adjoint is an algebra homomorphism.}
\label{fig:mult}
\end{figure}
In the rest of this paper, equality using the braided commutativity is denoted by $\underset{\text{(bc)}}{=}$.  
\end{proof}
\begin{proposition}
The adjoint coaction commutes with the braiding $\Psi$ as follows.  
\begin{equation}
(id\otimes id \otimes \mu)(id \otimes \Psi \otimes id)(\ad \otimes \ad) \, \Psi
=
(\Psi \otimes id)(id\otimes id \otimes \mu)(id \otimes \Psi \otimes id)(\ad \otimes \ad).
\end{equation}
\label{prop:adproduct}
\end{proposition}
\begin{proof}
This relation comes from the braided commutativity as explained in Figure~\ref{fig:adbraiding}.  
\begin{figure}[htb!]
\[
\includegraphics[scale=0.7]{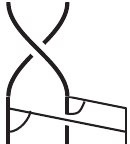}
\ \ \raisebox{10mm}{=}\ \ 
\includegraphics[scale=0.7]{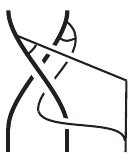}
\ \ 
\raisebox{10mm}{$\underset{\text{(bc)}}{=}$}
\ \ 
\includegraphics[scale=0.7]{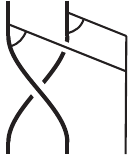}
\]
\caption{The adjoint coaction is commutative with the braiding.}
\label{fig:adbraiding}
\end{figure}
\end{proof}
\begin{proposition}
The adjoint coaction commutes with the antipode $S$, i.e.
\begin{equation}
\ad \circ S = (S \otimes id)\circ \ad.  
\label{eq:adS}
\end{equation}
\label{prop:adS}
\end{proposition}
\begin{proof}
This relation comes from the braided comutativity as explained in Figure~\ref{fig:adS}. 
The braided commutativity is used at the second equality. In the fourth equality we used the antipode axiom and
in the final equation the axiom relating $S$ and multiplication as used.  
\end{proof}
\begin{proposition}
$\mu \circ (id \otimes S) \circ \ad = S^2$.  
\label{prop:adS1}
\end{proposition}
\begin{proof}
This comes from the equalities of diagrams in Figure \ref{fig:adSS}. 
\end{proof}
\section{Representation of braid groups}
In this section, we recall the representation of
the braid group $B_n$  to $\End(A^{\otimes n})$ constructed by using the adjoint action of $A$.  
To construct  representations of braid groups, $A$ is not required to be braided commutative.  
However, for the distributivity of the representation given by Proposition \ref{prop:product}, $A$ has to be braided commutative.  
\subsection{Representation of generators}
The braid group $B_n$ is defined by the following generators and relations.  
\begin{multline}
B_n \!= 
\big\langle
\sigma_1, \sigma_2, \cdots, \sigma_{n-1} \mid
\sigma_i\,\sigma_{i+1}\, \sigma_i = 
\sigma_{i+1} \, \sigma_i \, \sigma_{i+1}, \ 
\sigma_i \, \sigma_j = \sigma_j \, \sigma_i
\\ 
(|i-j| \geq 2)
\big\rangle.
\label{eq:braid}
\end{multline}
\par
We define a braided Hopf diagram corresponding to the braid generators by generalizing the definition of $R_{\rm ad}$ in \cite{CS}, which is based on \cite{W2}.  
These are braided version of the Wirtinger presentation for the fundamental group of a knot complement.  
\begin{figure}[htb!]
\[
\begin{matrix}
R =
\begin{matrix}
\includegraphics[scale=0.6]{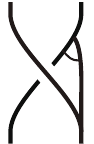}
\end{matrix}\ ,
& 
\qquad\qquad
&
R^{-1}  = \ 
\begin{matrix}
\includegraphics[scale=0.6]{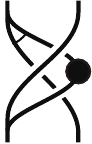}
\end{matrix}
\end{matrix}
\]
\caption{$R$, $R^{-1}$ for crossings}
\label{fig:R}
\end{figure}
For $\sigma_i ^{\pm1}\in B_n$,  let 
\begin{equation}
\rho(\sigma_i^{\pm1}) = id^{\otimes (i-1)} \otimes R^{\pm1} \otimes id^{\otimes(n-i-1)} \in \End(A^{\otimes n}),
\label{eq:R}
\end{equation} 
where $R^{\pm1} : A^{\otimes2} \to A^{\otimes2}$ is given in Figure \ref{fig:R}.   
\par
\begin{figure}[htb!]
\[
\begin{matrix}
&
\boldsymbol x \quad \otimes \quad \boldsymbol y
& & &\boldsymbol x
\\
\boldsymbol{\mu}^{(m)} \ : \ &
\begin{matrix}
\includegraphics[scale=0.7]
{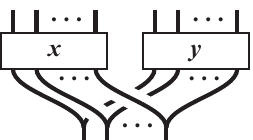}
\end{matrix}\ , &\qquad &
\boldsymbol{\Delta}^{(m)} \ : &
\begin{matrix}
\includegraphics[scale=0.7]
{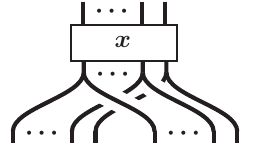}
\end{matrix} 
\\
&\boldsymbol{\mu}^{(k)}(\boldsymbol x \otimes \boldsymbol y) & & &
\boldsymbol{\Delta}^{(m)}(\boldsymbol{x})
\end{matrix}
\]
\caption{The generalized multiplication $\boldsymbol{\mu}^{(m)} : A^{\otimes m} \otimes A^{\otimes m} \to A^{\otimes m}$ and the generalized coproduct $\boldsymbol{\Delta}^{(m)} :  A^{\otimes m} \to A^{\otimes m} \otimes A^{\otimes m}$.}
\label{fig:genmult}
\end{figure}

In the rest, we  use the following operators 
$\Delta_i$, $\mu_i$ $S_i$, $\Psi_i$, $\varepsilon_i$ acting on $A^{\otimes n}$. 
They are given by the following.    
\begin{align*}
\Delta_i &= id^{i-1}\otimes \Delta \otimes id^{n-i}, 
\qquad
\mu_i = id^{i-1}\otimes \mu \otimes id^{n-i-1}, 
\\ 
S_i &=id^{i-1}\otimes S \otimes id^{n-i}, 
\qquad
\Psi_i  =   id^{i-1}\otimes \Psi \otimes id^{n-i-1}
\\
\varepsilon_i &= 
id^{i-1}\otimes \varepsilon \otimes id^{n-i}.
\end{align*}
%
We also use the generalized multiplication $\boldsymbol{\mu}^{(m)} : A^{\otimes m} \otimes A^{\otimes m} \to A^{\otimes m}$ and the generalized coproduct $\boldsymbol{\Delta}^{(m)} :  A^{\otimes m} \to A^{\otimes m} \otimes A^{\otimes m}$ given by
the  diagrams in Figure \ref{fig:genmult}.  
\subsection{Adjoint coaction}
\begin{figure}[htb!]
\[
\begin{matrix}
\text{\footnotesize $x_1\, x_2 \quad x_n$\qquad}
\\
\includegraphics[scale=0.7]{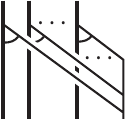}
\end{matrix}
\]
\caption{Adjoint coaction $\Ad: A^{\otimes n} \to A^{\otimes n}\otimes A$.}
\label{fig:Ad}
\end{figure}

We define an adjoint coaction 
$
\Ad : A^{\otimes n} \to A^{\otimes n} \otimes A
$
as in Figure~\ref{fig:Ad}.
\begin{proposition}
The adjoint coaction $\Ad$ commutes with $\rho(b)$ for $b \in B_n$, i.e.
\begin{equation}
\Ad \circ \rho(b) = (\rho(b) \otimes id) \circ \Ad.
\label{eq:adbraid}
\end{equation}
\label{prop:adbraid}
\end{proposition} 
\begin{proof}
This comes from the following commutativity of $R$ and $\Ad$.  
\begin{equation}
\Ad \circ R = (R \otimes id) \circ \Ad.
\label{eq:adR}
\end{equation}
This is proved by the graphical computation in Figure~\ref{fig:Adproof}.  
\begin{figure}[htb!]
\[
\begin{matrix}
\includegraphics[scale=0.5]{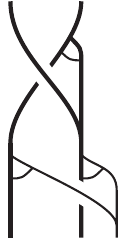}
\end{matrix}
\ =\ 
\begin{matrix}
\includegraphics[scale=0.5]{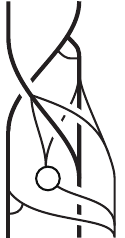}
\end{matrix}
\ =\ 
\begin{matrix}
\includegraphics[scale=0.5]{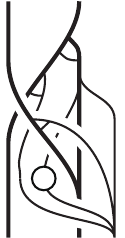}
\end{matrix}
\ =\ 
\begin{matrix}
\includegraphics[scale=0.5]{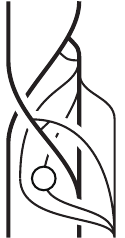}
\end{matrix}
\ =\ 
\begin{matrix}
\includegraphics[scale=0.5]{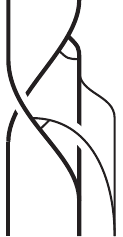}
\end{matrix}
\ =\ 
\begin{matrix}
\includegraphics[scale=0.5]{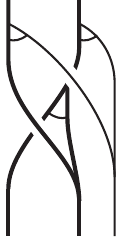}
\end{matrix}
\]
\caption{Commutativity of $\Ad$ and $R$. } 
\label{fig:Adproof}
\end{figure}
\end{proof}
\subsection{Representation of braid groups}
Now we construct a representation of braid groups in $\End(A^{\otimes n})$ 
\begin{theorem}[\cite{W2}, Proposition 1]\label{th:braid}
The map $\rho$ defined for generators of $B_n$ in \eqref{eq:R} extends to an algebra homomorphism from the group algebra ${\mathbb C}B_n$ to $\End(A^{\otimes n})$.  
\end{theorem}
\begin{proof}
We first show that $\rho(\sigma_i) \, \rho(\sigma_i^{-1}) = \rho(\sigma_i^{-1}) \, \rho(\sigma_i) = 1$.  
To show these, we prove
$R \, R^{-1} = R^{-1} \, R = id \otimes id$ by the graphical computation in Figure~\ref{fig:RRm} and Figure~\ref{fig:RmR}.  
\begin{figure}[htb!]
\[
\begin{matrix}
R^{-1}
\\[30pt]
R
\end{matrix}\quad
\begin{matrix}
\includegraphics[scale=0.45]{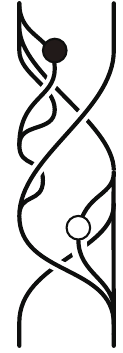}
\end{matrix}
\ =\ 
\begin{matrix}
\includegraphics[scale=0.45]{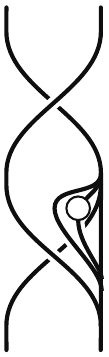}
\end{matrix}
\ =\ 
\begin{matrix}
\includegraphics[scale=0.45]{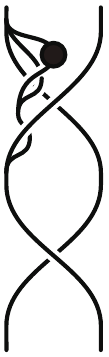}
\end{matrix}
\ =\ 
\begin{matrix}
\includegraphics[scale=0.45]{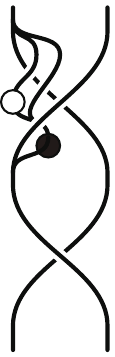}
\end{matrix}
\ =\ 
\begin{matrix}
\includegraphics[scale=0.45]{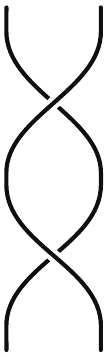}
\end{matrix}
\ =\ 
\begin{matrix}
\includegraphics[scale=0.45]{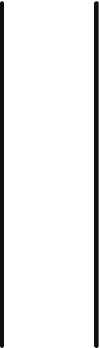}
\end{matrix}
\]
\caption{$R \, R^{-1} =  id \otimes id$.}
\label{fig:RRm}
\end{figure}
\begin{figure}[htb!]
\[
\begin{matrix}
R
\\[30pt]
R^{-1}
\end{matrix}
\begin{matrix}
\includegraphics[scale=0.45]{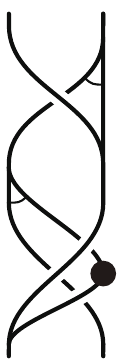}
\end{matrix}
\ =\ 
\begin{matrix}
\includegraphics[scale=0.45]{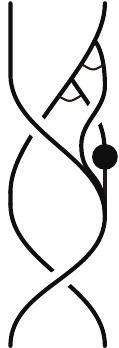}
\end{matrix}
\ =\ 
\begin{matrix}
\includegraphics[scale=0.45]{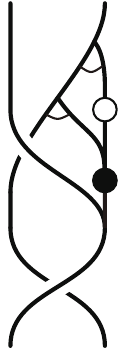}
\end{matrix}
\ =\ 
\begin{matrix}
\includegraphics[scale=0.45]{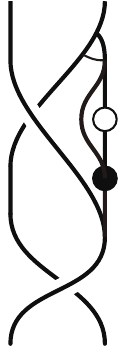}
\end{matrix}
\ =\ 
\begin{matrix}
\includegraphics[scale=0.45]{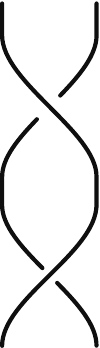}
\end{matrix}
\ =\ 
\begin{matrix}
\includegraphics[scale=0.45]{RRm5}
\end{matrix}
\]
\caption{$R^{-1} \, R =  id \otimes id$.}
\label{fig:RmR}
\end{figure}
The braid relation 
$\sigma_i \, \sigma_{i+1} \, \sigma_i 
= 
\sigma_{i+1} \, \sigma_i \, \sigma_{i+1}$ 
comes from 
\[
(R \otimes id)(id \otimes R)(R \otimes id)
=
(id \otimes R)(R \otimes id)(id \otimes R),
\]
which is shown by the graphical computation in Figure~\ref{fig:RIII}.  
\begin{figure}[htb!]
\[
\begin{matrix}
\includegraphics[scale=0.5]{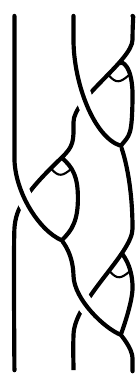}
\end{matrix}
\ =\ 
\begin{matrix}
\includegraphics[scale=0.5]{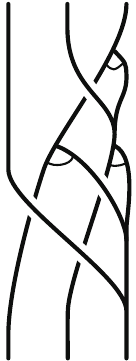}
\end{matrix}
\ \underset{\eqref{eq:adR}}{=}\ 
\begin{matrix}
\includegraphics[scale=0.5]{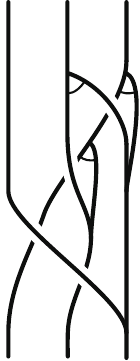}
\end{matrix}
\ =\ 
\begin{matrix}
\includegraphics[scale=0.5]{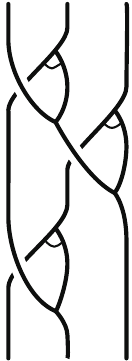}
\end{matrix}
\]
\caption{Braid relation $(R \otimes id)(id \otimes R)(R \otimes id)
=
(id \otimes R)(R \otimes id)(id \otimes R)$.}
\label{fig:RIII}
\end{figure}
We also have $\sigma_i\,\sigma_j = \sigma_j\,\sigma_i$ for $j-i \geq 2$ since
\[
R_i \, R_j
=
id ^{\otimes(i-1)} \otimes R \otimes id^{\otimes(n-i-j-2)} \otimes R \otimes id^{\otimes(n-j-1)}
=
R_j \, R_i.
\]
where $R_i = id ^{\otimes(i-1)} \otimes R \otimes id^{\otimes(n-i-1)}$.  
Hence the relations of $B_n$ are all satisfied.  
\end{proof}
%
%
%
\subsection{Distributivity of $\rho(b)$.}
\begin{figure}[htb]
\[
\begin{matrix}
\includegraphics[scale=0.7]{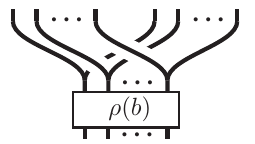}
\end{matrix}
\ = \ 
\begin{matrix}
\includegraphics[scale=0.7]{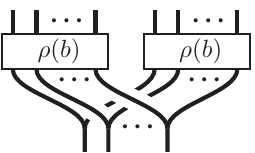}
\end{matrix}
\]
\caption{The representation $\rho$ is distributive.}
\label{fig:distributive}
\end{figure}

\begin{figure}[htb]
\[
\begin{matrix}
\includegraphics[scale=0.6]{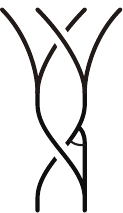}
\end{matrix}
\underset{\text{Prop. \ref{prop:admult}}}{=}
\begin{matrix}
\includegraphics[scale=0.6]{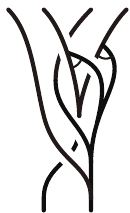}
\end{matrix}
\underset{\text{(bc)}}{=}
\begin{matrix} 
\includegraphics[scale=0.6]{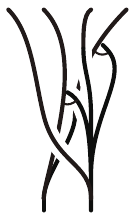}
\end{matrix}
=
\begin{matrix} 
\includegraphics[scale=0.6]{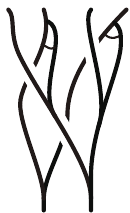}
\end{matrix}
\]
\caption{$R$ is distributive over the multiplication.}
\label{fig:mu}
\end{figure}

The representation $\rho(b)$ is distributive over the multiplication as follows.  
\begin{proposition}
Assume that $A$ is braided commutative.
For $\boldsymbol{x}, \boldsymbol{y} \in A^{\otimes n}$,  we have
\begin{equation}
\rho(b) \, \boldsymbol{\mu}^{(n)}(\boldsymbol x \otimes \boldsymbol y)
=
\boldsymbol{\mu}^{(n)}\big(\rho(b) \, \boldsymbol x\otimes \rho(b) \, \boldsymbol y\big).
\label{eq:split}
\end{equation}
This relation is explained graphically in Figure \ref{fig:distributive}.
\label{prop:product}
\end{proposition}
%
\begin{proof}
It is enough to show that
\[
R \, \boldsymbol{\mu}^{(2)}(\boldsymbol x \otimes \boldsymbol y) 
= 
\boldsymbol{\mu}^{(2)} \, (R \otimes R)(\boldsymbol x \otimes \boldsymbol y)
\]
for the multiplication $\boldsymbol{\mu}^{(2)} : A^{\otimes 2} \otimes A^{\otimes 2} \to A^{\otimes 2}$ and $\boldsymbol{x} = x_1 \otimes x_2,\ \boldsymbol{y}=y_1 \otimes y_2\in A^{\otimes 2}$, 
which is proved graphically in Figure~\ref{fig:mu}.  
\end{proof}
\section{Space of braided Hopf algebra representations of a knot}
Throughout this section $A$ is a braided Hopf algebra that is braided commutative.  
For any knot $K$, we construct the space of $A$ representations of $K$ as a quotient of $A^{\otimes n}$ by a module $I_b$ determined by a braid $b\in B_n$ whose closure is $K$.  
The number $n$ and the module $I_b$ depend on the choice of the braid $b$, but it is shown that the resulting quotient $A^{\otimes n}/I_b$ are isomorphic if the closures of the braids are isotopic.  Moreover they are isomorphic as $Ad$ comodules.
\subsection{Knots as braid closures}
\begin{figure}[htb]
\[
\begin{matrix}
\begin{matrix}\includegraphics[scale=0.65]{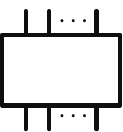}
\hspace{-8mm}\raisebox{5.5mm}{$b$}\hspace{6mm}\ 
\end{matrix}
&\ \  \longrightarrow \  \ 
&
\begin{matrix}
\begin{matrix}\includegraphics[scale=0.65]{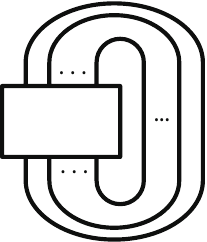}
\hspace{-17mm}\raisebox{12mm}{$b$}\hspace{17mm}\ 
\end{matrix}
\end{matrix}
\\
b & &  \text{$\widehat b$ : the closure of $b$}
\end{matrix}
\]
\caption{ Braid closure.}
\label{fig:closure}
\end{figure}

Let $K$ be a knot in $S^3$, then it is known that there is a braid $b\in B_n$ for certain $n\in \mathbb{N}$ such that $K$ is isotopic to the closure of $b$.  
The closure of $b$ is obtained by connecting the top points and bottom points of $b$ as in Figure~\ref{fig:closure}, and is denoted by $\widehat b$.  
\subsection{Space of $A$ representations}
\begin{figure}[htb]
\[
\begin{matrix}
d({b}) : \ 
\begin{matrix}
\includegraphics[scale=0.7]{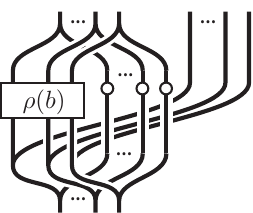}
\end{matrix}\ ,
& & 
\varepsilon^{\otimes n}\otimes id^{\otimes n} : \ 
\begin{matrix}
\includegraphics[scale=0.7]{braidclosure1}
\end{matrix}
\end{matrix}
\]
\caption{The braided Hopf diagrams for  $d({b})$ and $\boldsymbol{\varepsilon}_n \otimes id^{\otimes n}$. }
\label{fig:braidHopfdiagram}
\end{figure}

For $b \in B_n$, let $d({b})$ be the  braided Hopf diagram given in Figure \ref{fig:braidHopfdiagram}.   
Then $d(b)$ is an element in $\Hom(A^{\otimes 2n}, A^{\otimes n})$.  
Let us assign $x_1$, $x_2$, $\cdots$, $x_n$, $y_1$ $\cdots$, $y_n \in A$ to the top points of $d({b})$, let $\boldsymbol{x} = x_1 \otimes \cdots\otimes x_n$, $\boldsymbol{y} = y_1 \otimes \cdots \otimes y_n$ and 
 $d({{b}})(\boldsymbol{x}\otimes \boldsymbol{y})$ be the image of $\boldsymbol{x}\otimes \boldsymbol{y}$ by $d(b)$, which is an element in $A^{\otimes n}$.  
Let 
\[
I_{d(b)} 
= 
\image\big(d({b}) - {\varepsilon}^{\otimes n}\otimes id^{\otimes n}\big).
\]
%
\begin{definition}
A submodule $I$ of $A^{\otimes n}$ is called an {\it $\Ad$-comodule} if $\Ad(I) \subset I \otimes A$.  A morphism between two $\Ad$-comodules $I$ and $J$ is a module map $f:I\to J$ that commutes with $\Ad$ in the sense that $\Ad \circ f = (f\otimes id_A) \circ \Ad$. 
\end{definition}
\begin{proposition}
$I_{d(b)}$ is an $\Ad$-comodule of $A^{\otimes n}$. 
\end{proposition}
\begin{proof}
From \eqref{eq:adcomult}, \eqref{eq:admult} and \eqref{eq:adbraid}, we have $\Ad\circ d(b) = (d(b) \otimes id) \circ \Ad$ and
$\Ad \circ ({\varepsilon}^{\otimes n}\otimes id^{\otimes n} )
=
({\varepsilon}^{\otimes n}\otimes id^{\otimes (n+1)})\circ \Ad$. 
Therefore, $\Ad \circ \big(d({b}) - {\varepsilon}^{\otimes n}\otimes id^{\otimes n}\big) 
= 
\big(d({b}) - {\varepsilon}^{\otimes n}\otimes id^{\otimes (n+1)}\big)\circ \Ad$ and the image of $\Ad \circ \big(d({b}) - {\varepsilon}^{\otimes n}\otimes id^{\otimes n}\big)$ is contained in $I_{d(b)} \otimes A$.  
\end{proof}
\par
Let $A_b = A^{\otimes n} / I_{d(b)}$, 
then $A_{b}$ is an $\Ad$-comodule of $A^{\otimes n}$ and it satisfies the following.   
\begin{theorem}
If the closures of two braids $b_1 \in B_{n_1}$ and $b_2 \in B_{n_2}$ are isotopic, then $A_{b_1}$ and $A_{b_2}$ are isomorphic $\Ad$-comodules. 
In other words, $A_b$ is an invariant of the knot (or link) $\widehat b$.  
\label{th:main}
\end{theorem}

\begin{definition}
The $\Ad$-comodule $A_{b}$ is called {\it the space of $A$ representations} of the closure $\widehat b$.  
\label{def:space}
\end{definition}
The $\Ad$-comodule structure on $A_b$ allows us to pass to the coinvariants. This should generalize the conjugation invariant functions on the
representation variety and hence we introduce the following definition.
\begin{definition}
We say  the {\it quantum $A$ character variety} of $K$ is the module of coinvariants $A_b{}^{A}$ under the coaction of $\Ad$ on $A_b$.  
\end{definition}
It should be noted that our quantum $A$ character variety is not an algebra but only a module. 

\subsection{Equivalent pairs} 
To prepare our proof of the main theorem, Theorem \ref{th:main}, we introduce the notion of equivalent pairs of Hopf diagrams.
\begin{definition}
Hopf diagrams $d_1$, $d_2 \in \BHD(m, n)$ are called  {\it  $b$ equivalent }  if $A^{\otimes n}/\big(d_1 -d_2)(A^{\otimes m}\big)$ 
are isomorphic to $A_b$ as
$\Ad$-comodules. When $d_1$ and $d_2$ are equivalent we will denote this by 
$d_1 \sim_b d_2$ .  
Especially, $d(b)\sim_b {\varepsilon}^{\otimes n}\otimes id^{\otimes n}$.  
\end{definition}
{
Let $M_n$ be the kernel of 
${\varepsilon}^{\otimes n}$ in $A^{\otimes n}$. For $d \in \mathrm{BHD}(2n, n)$ we denote the induced map from $\big(A^{\otimes n}/M_n\big) \otimes A^{\otimes n}$ to $A^{\otimes n}/d(M_n \otimes A^{\otimes n})$ by $\bar d$.  
\begin{lemma}
If $d \in \mathrm{BHD}(2n,n)$ satisfies $d(1^{\otimes n}\otimes \by) = \by$, 
then
the image of 
$d - {\varepsilon}^{\otimes n}\otimes id^{\otimes n}$ 
is equal to $d(M_n\otimes A^{\otimes n})$.  
Especially,  
$I_{d(b)}= d(b)(M_n\otimes A^{\otimes n})$.  
\label{lem:image}
\end{lemma} 
\begin{proof}
The assumption $d(1^{\otimes n}\otimes \by) = \by$  
implies that $\big(d-{\varepsilon}^{\otimes n}\otimes id^{\otimes n}\big)(1^{\otimes n}\otimes \by) = 0$.  
Since 
$A^{\otimes n} = \mathbb{C}\big(1^{\otimes n} \big)\oplus M_n$  
and 
${\varepsilon}^{\otimes n}(M_n) = 0$, 
we get
$\big(d- {\varepsilon}^{\otimes n}\otimes id^{\otimes n}\big)(A^{\otimes n}\otimes A^{\otimes n}) 
= 
\big(d-{\varepsilon}^{\otimes n}\otimes id^{\otimes n}\big)(M_n\otimes A^{\otimes n}) 
=  
d(M_n\otimes A^{\otimes n})$.   
The last statement comes from the fact that $d(b)(1^{\otimes n}\otimes \boldsymbol{y}) = \boldsymbol{y}$.  
\end{proof}
}
%
%
\subsection{Moves at the top of BHD}
Let $d$ be an element of  $\BHD(2n,n)$ which is $b$ equivalent to $\varepsilon^{\otimes n} \otimes id^{\otimes n}$.  
The following proposition shows that we can modify 
$d$ at 
certain kinds of  multiplications, adjoints and braidings near the top of the diagram $d$ so that the corresponding $\overline{d}$ does not change. 
\begin{proposition}
Let  $d_1$, $d_2 \in \BHD(2n, n)$ be a pair of braided Hopf diagrams shown below.  
Assuming that $d_j(1^{\otimes n} \otimes \by) = \by$ for $j=1$, $2$, we have
 $d_1 \sim_b {\varepsilon}^{\otimes n}\otimes id^{\otimes n}$ if and only if $d_2 \sim_b {\varepsilon}^{\otimes n}\otimes id^{\otimes n}$.  
In the pictures the index $i$ is in $\{1,2,\cdots, n\}$.
\label{prop:BHDmult}
\begin{align*}
\text{\bf HmL}&\qquad
d_1 \ \begin{matrix}
\includegraphics[scale=0.6]{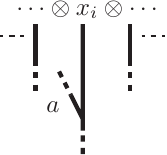}
\end{matrix}
\ \ \longleftrightarrow \quad  d_2 \ 
\begin{matrix}
\includegraphics[scale=0.6]{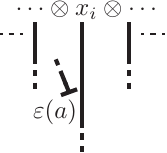}
\end{matrix}
\\
\text{\bf HmR}&
\qquad
d_1 \ 
\begin{matrix}
\includegraphics[scale=0.6]{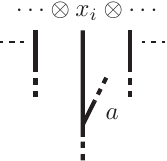}
\end{matrix}
\ \ \longleftrightarrow \quad  d_2 \ 
\begin{matrix}
\includegraphics[scale=0.6]{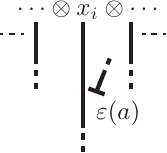}
\end{matrix}
\\
\text{\bf Ha}&\qquad
d_1 \ 
\begin{matrix}
\includegraphics[scale=0.6]{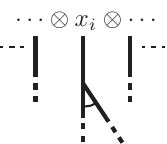}
\end{matrix}
\ \  \longleftrightarrow \quad   d_2 \ 
\begin{matrix}
\includegraphics[scale=0.6]{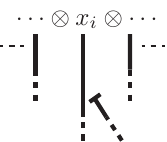}
\end{matrix}
\\
\text{\bf HmL'}&\qquad
d_1 \ \begin{matrix}
\includegraphics[scale=0.6]{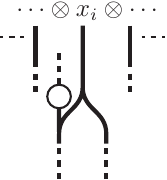}
\end{matrix}
\ \ \longleftrightarrow \quad  d_2 \ 
\begin{matrix}
\includegraphics[scale=0.6]{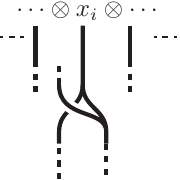}
\end{matrix}
\\
\text{\bf HmR'}&
\qquad
d_1 \ 
\begin{matrix}
\includegraphics[scale=0.6]{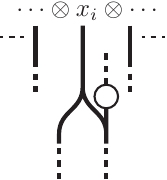}
\end{matrix}
\ \ \longleftrightarrow \quad d_2 \ 
\begin{matrix}
\includegraphics[scale=0.6]{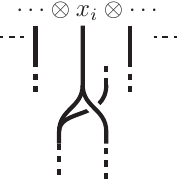}
\end{matrix}
\\
\text{\bf Ha'}&\qquad
d_1 \ 
\begin{matrix}
\includegraphics[scale=0.6]{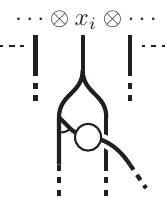}
\end{matrix}
\ \ \ \longleftrightarrow \quad  d_2 \ 
\begin{matrix}
\includegraphics[scale=0.6]{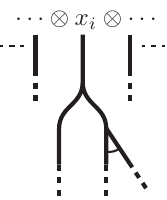}
\end{matrix}
\\
\text{\bf HcL}&\qquad
d_1 \quad
\begin{matrix}
\includegraphics[scale=0.6]{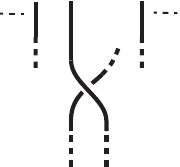}
\end{matrix}
\ \ \longleftrightarrow \quad 
d_2 \ \  
\begin{matrix}
\includegraphics[scale=0.6]{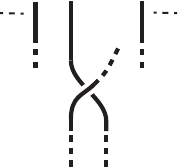}
\end{matrix}
\\
\text{\bf HcR}&\qquad
d_1 \quad 
\begin{matrix}
\includegraphics[scale=0.6]{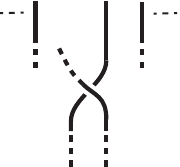}
\end{matrix}
\ \ \longleftrightarrow \quad\ d_2 \ \  
\begin{matrix}
\includegraphics[scale=0.6]{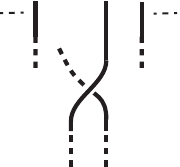}
\end{matrix}
\\
\text{\bf Hf}&\qquad
d_1 \quad 
\begin{matrix}
\includegraphics[scale=0.6]{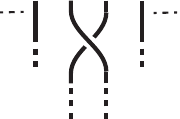}
\end{matrix}
\ \ \longleftrightarrow \quad\ d_2 \ \  
\begin{matrix}
\includegraphics[scale=0.6]{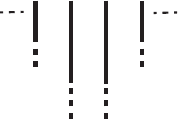}
\end{matrix}
\end{align*}
\end{proposition} 
%
\begin{proof}
{\bf HmL, HmR:}  
Let $d_1$, $d_2$ be the braided Hopf diagrams of {\bf HmL},
 $I_1=d_1(M_n\otimes A^{\otimes n})$ and $I_2= d_2(M_n\otimes A^{\otimes n})$.  
Let $\overline{d_1}$ be the map from $A^{\otimes 2n}/\ker \varepsilon_i$ to $A^{\otimes n}/I_1$ induced by $d_1$. 
The multiplication of $x_i$ by $a$ in the diagram $d_1$ is the multiplication of $\overline{x}_i\in A/M_1$ by $a$, which is equal to $\varepsilon(a) \, \overline{x}_i$.  
Therefore, $\overline{d}_1 = \overline{d}_2$, where $\overline{d}_2$ is the map from $A^{\otimes 2n}/\ker \varepsilon_i$ to $A^{\otimes n}/I_1$  given by $d_2$.  
This implies 
that $I_2 \subset I_1$.  

On the other hand, 
let ${\widetilde{d_2}}$ be the map from $A^{\otimes 2n}/\ker \varepsilon_i$ to $A^{\otimes n}/I_2$ induced by $d_2$. 
The multiplication of $a$ and $x_i$ in the diagram $d_1$ is the multiplication of $\overline{x}_i\in A/M_1$ by $a$, which is equal to $\varepsilon(a) \, \overline{x}_i$.  
Therefore, $\widetilde{d}_2 = \widetilde{d}_1$, where $\widetilde{d}_1$ is the map from $A^{\otimes 2n}/\ker \varepsilon_i$ to $A^{\otimes n}/I_2$  given by $d_1$.  
This implies 
 $I_1 \subset I_2$.  
Hence we get $I_1 = I_2$.  
This means that $I_1 = I_{d(b)}$ if and only if $I_2 = I_{d(b)}$, which implies  $d_1 \sim_b {\varepsilon}^{\otimes n}\otimes id^{\otimes n}$ if and only if $d_2 \sim_b {\varepsilon}^{\otimes n}\otimes id^{\otimes n}$ by Lemma \ref{lem:image}.  
The proof for $d_1$, $d_2$ in {\bf HmR} is similar.  
\par
{\bf Ha:}  
Let $d_1$, $d_2$ be the braided Hopf diagrams of {\bf Ha},
 $I_1=d_1(M_n\otimes A^{\otimes n})$ and $I_2= d_2(M_n\otimes A^{\otimes n})$.  

Since $(\varepsilon \otimes id)\ad(M_1) =  (1 \otimes 1)\varepsilon(M_1) = 0$, $\ad(M_1)$ is contained in $M_1 \otimes A$.  
Hence $\ad$ induces the map $\overline{\ad}$ from $A/M_1$ to $\big(A/M_1\big) \otimes A$.  
On the other hand, $A/M_1$ is spanned by $1$ and $\ad(1) = 1 \otimes 1$, $\overline{\ad}(\overline{x}) = \overline{x} \otimes 1$ for $\overline{x} \in A/M_1$.  
This relation means that $\overline{d}_1 = \overline{d_2}$ as a map from 
$A^{\otimes 2n}/\ker \varepsilon_i$ to $A^{\otimes n}/I_1$.  
Therefore, $I_2 \subset I_1$.  
By exchanging the role of $d_1$ and $d_2$ as in the case of {\bf HmL}, we get $I_1 \subset I_2$.  
Hence $I_1 = I_2$.  
Therefore, $d_1 \sim_b {\varepsilon}^{\otimes n}\otimes id^{\otimes n}$ if and only if $d_2 \sim_b {\varepsilon}^{\otimes n}\otimes id^{\otimes n}$.  

\par
{\bf HmL', HmR':}  
Let $d_1$, $d_2$ be the braided Hopf diagrams of {\bf HmL'}.  
Assume that $d_1 \sim_b {\varepsilon}^{\otimes n}\otimes id^{\otimes n}$.  
In the pictures the symbol $a \longleftrightarrow b$ means that
$a$ and $b$ are both $b$ equivalent to  ${\varepsilon}^{\otimes n}\otimes id^{\otimes n}$.
By {\bf HmL}, we have a sequence of equalities as in Figure \ref{fig:BHDswitchproof} as a map from $\big(A^{\otimes n}/M_n\big)\otimes A^{\otimes n}$ to $A_b$.  
Hence $d_1(M_m\otimes A^{\otimes n}) = d_2(M_n\otimes A^{\otimes n})$ and $d_2 \sim_b {\varepsilon}^{\otimes n}\otimes id^{\otimes n}$.  
The proof for {\bf HmR'} is similar.  
%
\begin{figure}[htb!]
\[
\begin{matrix}
\begin{matrix}
\includegraphics[scale=0.6]{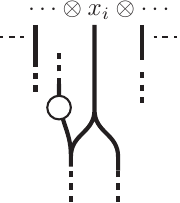}
\end{matrix}
\!=\!
\begin{matrix}
\includegraphics[scale=0.6]{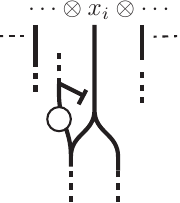}
\end{matrix}
\underset{\text{\bf HmL}}{\longleftrightarrow}
\begin{matrix}
\includegraphics[scale=0.6]{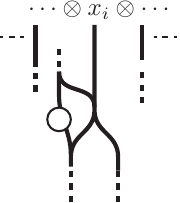}
\end{matrix}
\!=\!
\begin{matrix}
\includegraphics[scale=0.6]{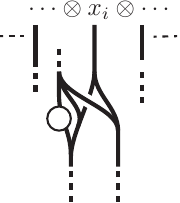}
\end{matrix}
\!=\!
\begin{matrix}
\includegraphics[scale=0.6]{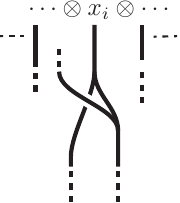}
\end{matrix}
\end{matrix}
\]
\caption{Proof for {\bf HmL'}.}
\label{fig:BHDswitchproof}
\end{figure}
%
%
\par
{\bf Ha':}  
Let $d_1$, $d_2$ be the braided Hopf diagrams of {\bf Ha'}.  
Assume that $d_1 \sim_b {\varepsilon}^{\otimes n}\otimes id^{\otimes n}$.  
By {\bf Ha}, we have a sequence of equalities in Figure \ref{fig:BHDadswitchproof} as a map from $\big(A^{\otimes n}/M_n\big)\otimes A^{\otimes n}$ to $A_b$.   In the third and fourth equalities we used Proposition \ref{prop.Adrels}.
Hence 
$d_1(M_m\otimes A^{\otimes n}) = d_2(M_n\otimes A^{\otimes n})$ 
and $d_2 \sim_b {\varepsilon}^{\otimes n}\otimes id^{\otimes n}$.  
The opposite direction is proved similarly.  
\begin{figure}[htb!]
\[
\begin{matrix}
\begin{matrix}
\includegraphics[scale=0.6]{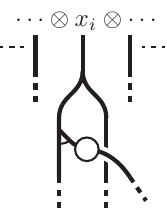}
\end{matrix}
=
\begin{matrix}
\includegraphics[scale=0.6]{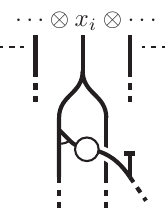}
\end{matrix}
\underset{\text{\bf Ha}}{\longleftrightarrow}
\begin{matrix}
\includegraphics[scale=0.6]{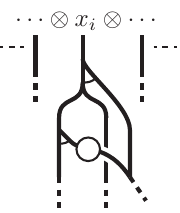}
\end{matrix}
=
\begin{matrix}
\includegraphics[scale=0.6]{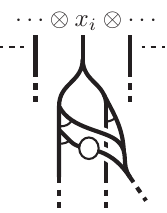}
\end{matrix}
=
\begin{matrix}
\includegraphics[scale=0.6]{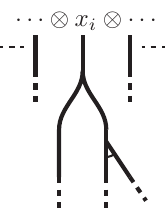}
\end{matrix}
\end{matrix}
\]
\caption{Proof for {\bf Ha'}.}
\label{fig:BHDadswitchproof}
\end{figure}
%
%
%
%
\par
{\bf HcL, HcR:}
Let $d_1$, $d_2$ be the braided Hopf diagrams of {\bf HcL}.  
Assume that $d_1 \sim_b {\varepsilon}^{\otimes n}\otimes id^{\otimes n}$.  
We have a sequence of equalities as in Figure \ref{fig:BHDPhiproof} as a map from $A^{\otimes n}/M_n\otimes A^{\otimes n}$ to $A_b$. 
In the rest of this paper, (bc)'s under equalities mean braided commutativity.  
Hence $d_1(M_n\otimes A^{\otimes n}) = d_2(M_n\otimes A^{\otimes n})$ and $d_2 \sim_b {\varepsilon}^{\otimes n}\otimes id^{\otimes n}$.  
The proof for {\bf HcR} is  given in Figure \ref{fig:BHDPhiproof2}.  
\begin{figure}[htb]
\[
\begin{matrix}
\begin{matrix}
\includegraphics[scale=0.6]{BHDPhi1}
\end{matrix}
\underset{\text{\bf Ha}}{\longleftrightarrow}
\begin{matrix}
\includegraphics[scale=0.6]{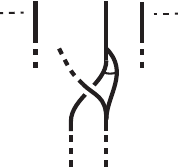}
\end{matrix}
\underset{\text{(bc)}}
{=}
\begin{matrix}
\includegraphics[scale=0.6]{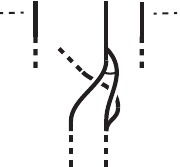}
\end{matrix}
\underset{\text{\bf Ha}}{\longleftrightarrow}
\begin{matrix}
\includegraphics[scale=0.6]{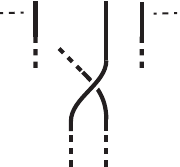}
\end{matrix}
\end{matrix}
\]
\caption{Proof for {\bf HcL}.}
\label{fig:BHDPhiproof}
\end{figure}
\begin{figure}[htb]
\[
\begin{matrix}
\includegraphics[scale=0.6]{BHDPhi4}
\end{matrix}
\underset{\text{\bf Ha}}
{\longleftrightarrow}
\begin{matrix}
\includegraphics[scale=0.6]{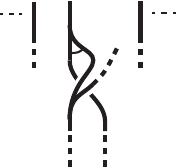}
\end{matrix}
\underset{\text{(bc)}}
{=}
\begin{matrix}
\includegraphics[scale=0.6]{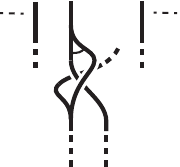}
\end{matrix}
\underset{\text{\bf Ha}}
{\longleftrightarrow}
\begin{matrix}
\includegraphics[scale=0.6]{BHDPhi3}
\end{matrix}
\]
\caption{Proof for {\bf HcR}.}
\label{fig:BHDPhiproof2}
\end{figure}
\par
{\bf Hf:} The relation {\bf Hf} comes from the facts that $\Psi$ is an automorphism of $A^{\otimes 2}$ and $\varepsilon^{\otimes 2} \circ \Psi = \varepsilon^{\otimes 2}$.  
\end{proof}
{
\begin{proposition}
Let $d$ be an element of  $\BHD(2n,n)$ which is $b$ equivalent to $\varepsilon^{\otimes n} \otimes id^{\otimes n}$.  
Let $d'$ be an element of  $\BHD(2n,2n)$ which gives an isomorphism from $A^{\otimes 2n}$ to $A^{\otimes 2n}$ such that $(\varepsilon^{\otimes n} \otimes id^{\otimes n})\circ d' = \varepsilon^{\otimes n} \otimes id^{\otimes n}$.  
Then $d \circ d'\sim_b\varepsilon^{\otimes n} \otimes id^{\otimes n}$. 
Especially,  for the diagram $d_{ij}$  having an arc connecting the  $(n+j)$-th strand to the $i$-th strand as in Figure \ref{fig:extramove}, then $d \circ d_{ij}\sim_b\varepsilon^{\otimes n} \otimes id^{\otimes n}$. 
Moreover, for the diagrams $d_1$ and $d_2$ in Figure \ref{fig:extramove}, $d_1 \sim_b {\varepsilon}^{\otimes n}\otimes id^{\otimes n}$ if and only if $d_2 \sim_b {\varepsilon}^{\otimes n}\otimes id^{\otimes n}$.  
\label{prop:extramove}
\end{proposition}
\begin{figure}[htb]
\[
\begin{matrix}
\text{\footnotesize\quad$i$\ \ \ \ \ \ \ \  $n+j$}
& & 
\text{\footnotesize\qquad$i$\ \ \ \ \ \  $n+j$}
& & 
\text{\footnotesize\qquad$i$\ \ \ \ \ \   $n+j$}
\\
\begin{matrix}
\includegraphics[scale=0.6]{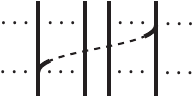}
\end{matrix} & &
\begin{matrix}
\includegraphics[scale=0.6]{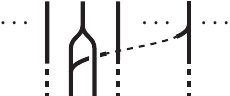} \end{matrix}
&\longleftrightarrow& 
\begin{matrix}
\includegraphics[scale=0.6]{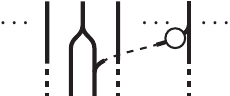}
\end{matrix}
\\
d_{ij}\in \BHD(2n, n) & & d_1 & & d_2
\end{matrix}
\]
\caption{The diagram $d_{ij}$ having an arc connecting the $(n+j)$-th strand to the $i$-th strand and the move comes from this diagram.}
\label{fig:extramove}
\end{figure}
\begin{proof}
Since $d'$ is an isomorphism, the image of 
$\big(d - \varepsilon^{\otimes n} \otimes id^{\otimes n}\big)\circ d' = d\circ d' - \varepsilon^{\otimes n} \otimes id^{\otimes n}$ is equal to $I_{d(b)}$.  
Therefore, $d\circ d'$ is b equivalent to $\varepsilon^{\otimes n} \otimes id^{\otimes n}$.  
The diagram $d_{ij}$ is an isomorphism since adding the antipode $S$ to the arc connecting the $(n+j)$-th strand to the $i$-th strand of $d_{ij}$, we get the inverse of $d_{ij}$.  
\par
By adding $d_{ij}$ to the top of $d_2$, we get $d_1$.  
This implies the last statement of the proposition.  
\end{proof}
}
\subsection{Another expression of $d(b)$}
Below we will show an alternative way of expressing the braided Hopf diagram $d(b)$ viewed as a map $\overline{d(b)}:A^{\otimes n}/M_n \otimes A^{\otimes n}\to A_b$.
Although this expression for $\overline{d(b)}$ will not be used in this text we include it because it corresponds more naturally to the closed braid of $b$. It also suggests that the construction given here for a closed braid may extend to a plat presentation of a knot. We will elaborate on this point in a future publication.
\begin{proposition}
For $b \in B_n$, the map $\overline{d(b)}$ induced by $d(b)$ satisfies
\begin{multline*}
\overline{d(b)} = 
(\mu_1\mu_3\cdots\mu_{2n-1} )\circ
\\
\big((\sigma_{2n-2}\sigma_{2n-4}\cdots\sigma_{2})(\sigma_{2n-3}\sigma_{2n-5}\cdots\sigma_{3})\cdots (\sigma_{n+1}\sigma_{n-1})\sigma_n\big)\circ
\\
(\boldsymbol{\mu}^{(n)}\otimes id^{\otimes n})\circ
(id^{\otimes n}\otimes \Psi^{(n)}) \circ
 (b \otimes S^{\otimes n}\otimes id^{\otimes n} )\\
\circ
\big(\sigma_n(\sigma_{n-1}^{-1}\sigma_{n+1}^{-1})\cdots (\sigma_3^{-1}\sigma_5^{-1}\cdots\sigma_{2n-3}^{-1}) (\sigma_2^{-1}\sigma_4^{-1}\cdots\sigma_{2n-2}^{-1})\big) \circ
\\
(\Delta_{2n-1}\Delta_{2n-3}\cdots \Delta_1)
\end{multline*}
as a map from $A^{\otimes n}/M_n\otimes A^{\otimes n}$ to $A_b$, where $\Psi^{(n)}$ is defined as
\[
\Psi^{(n)}= (\Psi_n \Psi_{n-1} \cdots \Psi_1)
(\Psi_{n+1} \Psi_{n} \cdots \Psi_2)\cdots
(\Psi_{2n-1} \Psi_{2n-2} \cdots \Psi_n).
\] 
$\Psi^{(n)}$ is a composition of $n^2$ maps $\Psi_j$, see also Figure \ref{fig:full}.
\label{prop:Pbraid}
\end{proposition}
\begin{proof}
We first remove the first and last $\sigma_{2n-2}$ as in Figure \ref{fig:removesigma}.
\begin{figure}[htb]
\begin{multline*}
\begin{matrix}
\includegraphics[scale=0.65]{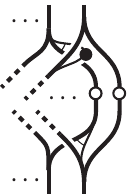}
\end{matrix}
=
\begin{matrix}
\includegraphics[scale=0.65]{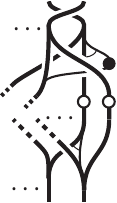}
\end{matrix}
\underset{\text{\bf Hf}}
{\ \longleftrightarrow\ }
\begin{matrix}
\includegraphics[scale=0.65]{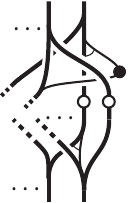}
\end{matrix}
\underset{\text{\bf HmR'}}
{\ \longleftrightarrow\ }
\begin{matrix}
\includegraphics[scale=0.65]{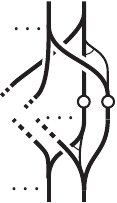}
\end{matrix}
=
\\
\begin{matrix}
\includegraphics[scale=0.65]{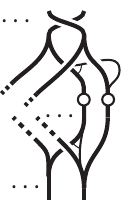}
\end{matrix}
\underset{\text{\bf Hf}}
{\ \longleftrightarrow\ }
\begin{matrix}
\includegraphics[scale=0.65]{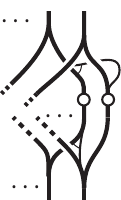}
\end{matrix}
\underset{\text{Prop. \ref{prop:adS}}}
{\ =\ }
\begin{matrix}
\includegraphics[scale=0.65]{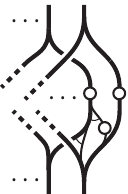}
\end{matrix}
{\ =\ }
\begin{matrix}
\includegraphics[scale=0.65]{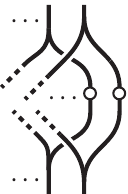}
\end{matrix}
\end{multline*}
\caption{Cancellation of the topmost $\sigma_{2n-2}^{-1}$ and the bottom-most $\sigma_{2n-2}$.}
\label{fig:removesigma}
\end{figure}
%
%
Then remove $\sigma_{2n-3}$, $\sigma_{2n-4}$, $\cdots$, $\sigma_{2n}$ similarly along the rightmost string.  
After doing these operations, do similar operations along the string next to the rightmost string.  
Repeat these operations for all $\sigma$'s along the strings with the antipode $S$. 
\end{proof}
\subsection{Markov moves}
It is  known that the closures of two braids $b_1 \in B_{n_1}$ and $b_2 \in B_{n_2}$ are isotopic in $S^3$ if and only if there is a sequence of the following two types of moves connecting $b_1$ to $b_2$.  
These moves are called the Markov moves and such $b_1$ and $b_2$ are called {\it Markov equivalent}.  
\begin{quote}
\item[\bf First Markov move (MI):]
\quad $b\, b' \longleftrightarrow b' \, b$ \quad for $b$, $b' \in B_n$.  
\item[\bf Second Markov move (MII):]
\quad $b \in B_n \longleftrightarrow \sigma_n^{\pm1} \, b   \in B_{n+1}$.  
\end{quote}
\begin{theorem}
 The quotient algebras $A_{b_1}$ and $A_{b_2}$ are isomorphic if $b_1$ and $b_2$ are Markov equivalent.  
 \end{theorem}
This theorem comes from Propositions \ref{prop:MI}, \ref{prop:MII1} and \ref{prop:MII2} in the following two subsections.  
%
%
%
 %
 \subsection{Invariance under the MI move}
First, we show that the quotient algebra keeps its structure when we apply an MI move.  
\begin{proposition}
For $b_1$, $b_2 \in B_n$, $A_{b_2b_1}$ is isomorphic to $A_{b_1b_2}$. 
\label{prop:MI} 
\end{proposition}
This comes from the following lemma.  
\begin{lemma}
For $b\in B_n$, $A_{\sigma_i b}$ is isomorphic to 
$A_{b \sigma_i}$. Also,
$A_{\sigma_i^{-1}b}$ is isomorphic to $A_{b\sigma_i^{-1}}$  
\end{lemma}
\begin{proof}
As a map from $A^{\otimes n}/M_n\otimes A^{\otimes n}$ to $A_{b\sigma_i}$, $d(\sigma_i b)$ is deformed as in Figures \ref{fig:MIdeform0} and \ref{fig:MIdeform} where $R$ and $R^{-1}$ are given in Figure \ref{fig:R}.
\begin{figure}[htb!]
\[
\begin{matrix}
\begin{matrix}
\epsfig{file=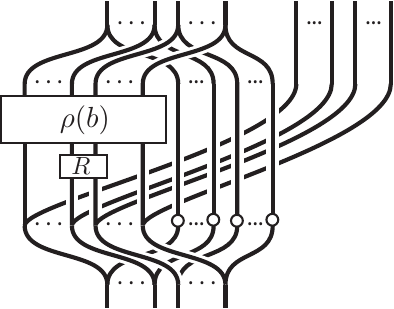, scale=0.6}
\end{matrix}
&\!\!\!\!\!\!
\underset{\text{Prop. \ref{prop:product}}}{=}
\!
&
\begin{matrix}
\epsfig{file=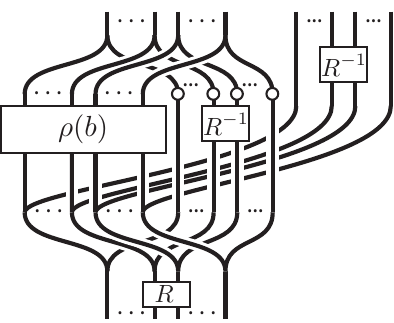, scale=0.6}
\end{matrix}
\\
d(\sigma_i b) & & d'
\end{matrix}
\]
\caption{Deform $d(b\sigma_i)$ to $d'$ as a map from $A^{\otimes n}/M_n$ to $A_{b\sigma_i}$.} 
\label{fig:MIdeform0}
\end{figure}
In Figure~\ref{fig:MIdeform}, $R^{-1}$ in Figure \ref{fig:MIdeform0} is moved upward and switched to the left strands by using the moves in Proposition \ref{prop:BHDmult}.  
After we apply the transformation suggested in Figure \ref{fig:MIdeform} to Figure \ref{fig:MIdeform0} we get Figure \ref{fig:MIdeform1}. It follows that
$\sigma_i\, d(b\sigma_i) (id^{\otimes n} \otimes \sigma_i^{-1})\sim_{\sigma_i b} {\varepsilon}^{\otimes n}\otimes id^{\otimes n}$.  
\begin{figure}[htb!]
\begin{multline*}
\begin{matrix}
\epsfig{file=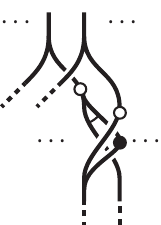, scale=0.53}
\end{matrix}
\underset{\text{Prop. \ref{prop:adS}}}{=}
\begin{matrix}
\epsfig{file=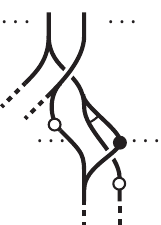, scale=0.53}
\end{matrix}
\underset{\text{\bf HcR}}{\longleftrightarrow}
\begin{matrix}
\epsfig{file=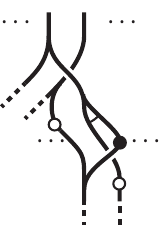, scale=0.53}
\end{matrix}
\underset{\text{\bf Ha'}}{\longleftrightarrow}
\begin{matrix}
\epsfig{file=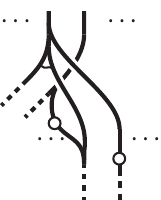, scale=0.53}
\end{matrix}
=
\begin{matrix}
\epsfig{file=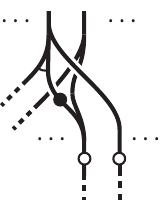, scale=0.53}
\end{matrix}
\underset{\text{\bf HmL'}}{\longleftrightarrow}
\\
\begin{matrix}
\epsfig{file=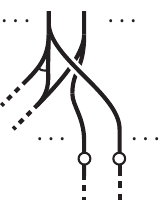, scale=0.53}
\end{matrix}
\underset{\text{(bc)}}{=}
\begin{matrix}
\epsfig{file=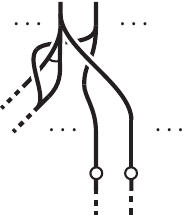, scale=0.53}
\end{matrix}
=
\begin{matrix}
\epsfig{file=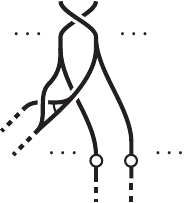, scale=0.53}
\end{matrix}
\underset{\text{\bf HcR}}{\longleftrightarrow}
\begin{matrix}
\epsfig{file=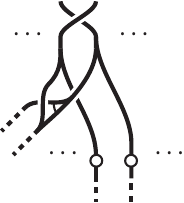, scale=0.53}
\end{matrix}
\underset{\text{\bf Hf}}{\longleftrightarrow}
\begin{matrix}
\epsfig{file=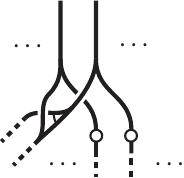, scale=0.53}
\end{matrix}
\end{multline*}
\hfill
\caption{Deform $d'$ as a map from $\big(A^{\otimes n}/M_n\big)\otimes A^{\otimes n}$ to $A_{b\sigma_i}$.}
\label{fig:MIdeform}
\end{figure}
Since 
$\sigma_i^{-1}$ 
is an automorphism of $A^{\otimes n}$, we have $\sigma_i \, d(b\sigma_i)\sim_{b\sigma_i}{\varepsilon}^{\otimes n}\otimes id^{\otimes n}$.  
This implies that  
\begin{multline*}
I_{d(\sigma_i b)}
=
 \big(\sigma_i\, d(b \sigma_i)  
 - 
  {\varepsilon}^{\otimes n}\otimes id^{\otimes n}\big)(A^{\otimes n}\otimes A^{\otimes n})\big\rangle 
  \\
=
 \big(\sigma_i\,d(b \sigma_i )   - 
 \sigma_i\,{\varepsilon}^{\otimes n}\otimes id^{\otimes n}\big)(A^{\otimes n}\otimes A^{\otimes n}) 
\\
=
 \sigma_i\,\big(d(b \sigma_i)  - 
 {\varepsilon}^{\otimes n}\otimes id^{\otimes n}\big)(A^{\otimes n}\otimes A^{\otimes n}) 
=
\sigma_i \, I_{d(b\sigma_i)}.
\end{multline*}
Hence the left multiplication of $\sigma_i^{-1}$ induces an isomorphism from $A^{\otimes n}/I_{d(\sigma_i b)}$ to $A^{\otimes n}/I_{d(b\sigma_i )}$. 
\par
For  $\sigma_i^{-1}b$, we  have
\[
A_{b\sigma_i^{-1}} = A_{\sigma_i\sigma_i^{-1}b\sigma_i^{-1}}\cong A_{\sigma_i^{-1}b\sigma_i^{-1}\sigma_i}= A_{\sigma_i^{-1}b}
\]
since $A_{b\sigma_i} \cong A_{\sigma_ib}$.  
\begin{figure}[htb!]
\[
\includegraphics[scale=0.6]{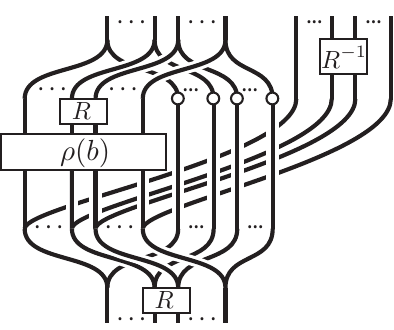}
\]
\caption{Result of the deformation of $d(\sigma_nb)$ in Figure \ref{fig:MIdeform}. }
\label{fig:MIdeform1}
\end{figure}
\end{proof}
\subsection{Invariance under  MII move}
\begin{figure}[htb!]
\[
\boldsymbol{\Psi}^{(n)} : 
\begin{matrix}
\includegraphics[scale=0.7]{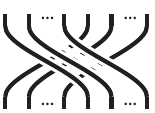}
\end{matrix}
\ ,
\qquad
\boldsymbol{\theta}^{(n)} : 
\begin{matrix}
\includegraphics[scale=0.7]{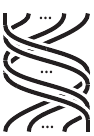}
\end{matrix} \ .
\]
\caption{$\boldsymbol{\Psi}^{(n)}$  and $\boldsymbol{\theta}^{(n)}$}
\label{fig:full}
\end{figure}

\begin{figure}[htb!]
\begin{multline*}
d(\sigma_n^{-1} b)
=\ 
\begin{matrix}
\includegraphics[scale=0.6]{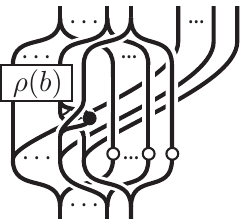}
\end{matrix}
=\ 
\begin{matrix}
\includegraphics[scale=0.65]{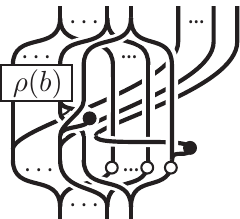}
\end{matrix}
\underset{\text{\bf Hf}}{\longleftrightarrow}
\ 
\begin{matrix}
\includegraphics[scale=0.65]{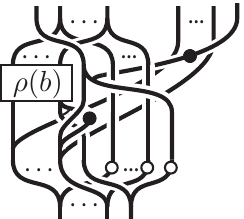}
\end{matrix}
\\
\underset{\text{Prop. \ref{prop:extramove}}}{\longleftrightarrow}
\ 
\begin{matrix}
\includegraphics[scale=0.65]{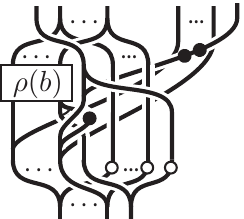}
\end{matrix}
\!\!\!\!
\underset{\text{\bf Hf}}{\longleftrightarrow}
\begin{matrix}
\includegraphics[scale=0.65]{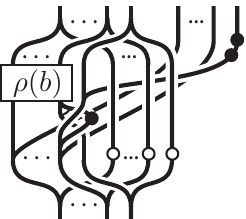}
\end{matrix}
\!\!\!\!
\underset{\text{(bc)}
}=\ 
\begin{matrix}
\includegraphics[scale=0.65]{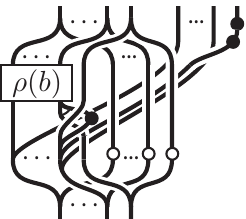}
\end{matrix}
\\
=
\begin{matrix}
\includegraphics[scale=0.65]{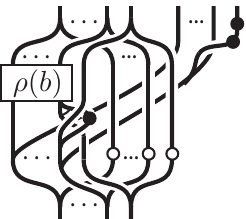}
\end{matrix}
=
d(\sigma_n^{-1}b)\otimes (id^{\otimes (n+1)}\otimes g\otimes 1)
\end{multline*}
\vspace{-5mm}
\caption{The map $d(\sigma_n^{-1} b)$ is equal to $d(\sigma_n^{-1}b)\otimes (id^{\otimes (n+1)}\otimes g\otimes 1)$.}
\label{fig:MII22}
\end{figure}

We show that the  quotient algebra keeps its structure by  MII move.   
We first compare $A_b$ and $A_{\sigma_n b}$.
\begin{proposition}
For $b\in B_n$, the $\Ad$-comocules  $A_{b}$ and $A_{\sigma_n b}$ are isomorphic.  
\label{prop:MII1}
\end{proposition}
\begin{proof}
Let $f$ be the linear surjection from $A^{\otimes (n+1)}$ to $A^{\otimes n}$ defined by 
\[
f = \mu_n\circ \Psi_n^{-1}.
\]
We first show that $f(I_{d(\sigma_n b)}) \subset I_{d(b)}$. 
From Figure \ref{fig:MII11},
\begin{figure}[htb!]
\begin{multline*}
f\circ d(\sigma_n b)
=
\begin{matrix}
\includegraphics[scale=0.65]{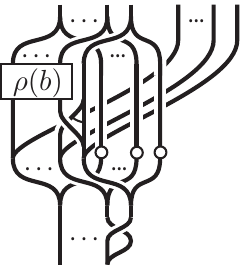}
\end{matrix}
\underset{\text{(bc)}}{=}
\begin{matrix}
\includegraphics[scale=0.65]{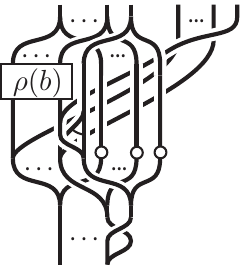}
\end{matrix}
=
\begin{matrix}
\includegraphics[scale=0.65]{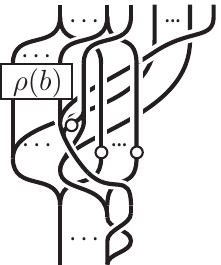}
\end{matrix}
\\
=
\begin{matrix}
\includegraphics[scale=0.65]{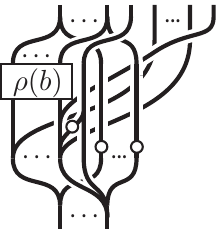}
\end{matrix}
=
\begin{matrix}
\includegraphics[scale=0.65]{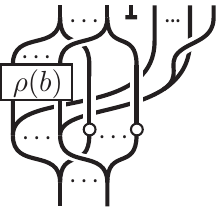}
\end{matrix}
=
d(b)\circ\varepsilon_{n+1}\circ(id^{\otimes (n+1)} \otimes f )
\end{multline*}
\vspace{-5mm}
\caption{The map $f\circ d(\sigma_n b)$ is equal to $d(b) \circ \varepsilon_{n+1}\circ(id^{\otimes (n+1)} \otimes f)$.}
\label{fig:MII11}
\end{figure}
 we know that
$
f\big(d(\sigma_n b)(\bx \otimes \by)\big) 
= 
d(b)\big(\varepsilon_{n+1}(\bx) \otimes f(\by)\big)
$,
and this means that $f(I_{d(\sigma_n b)}) \subset I_{d(b)}$.  
So $f$ induces a map $\overline{f}$ from $A_{d(\sigma_n b)}$ to $A_{d(b)}$.  
\par
Next, we show that $ \overline{f}$ is an isomorphism.  
Since $f$ is surjective, it is enough to show that $\overline{f}$ is injective.  
For this, we check that $\ker f$ is contained in $I_{d(\sigma_n b)}$.  
For $\bx \in A^{\otimes (n+1)}$ and $\by\in \ker f$,  
Figure \ref{fig:MII12} 
\begin{figure}[htb!]
\begin{multline*}
d(\sigma_n b)
=
\begin{matrix}
\includegraphics[scale=0.65]{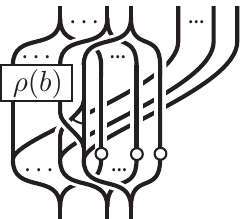}
\end{matrix}
\underset{\text{(bc)}}{=}
\begin{matrix}
\includegraphics[scale=0.65]{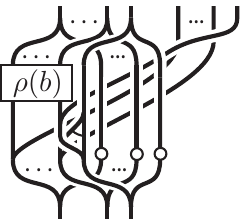}
\end{matrix}
=
\begin{matrix}
\includegraphics[scale=0.65]{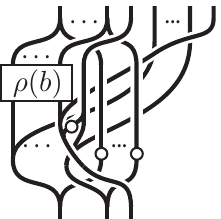}
\end{matrix}
\\
=
\begin{matrix}
\includegraphics[scale=0.65]{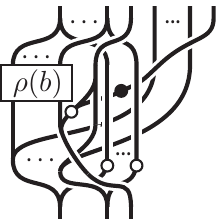}
\end{matrix}
\underset{\text{\bf Hf}}{\longleftrightarrow}
\begin{matrix}
\includegraphics[scale=0.65]{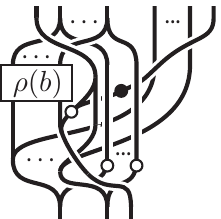}
\end{matrix}
\underset{\text{Prop. \ref{prop:extramove}}}{\longleftrightarrow}
\begin{matrix}
\includegraphics[scale=0.65]{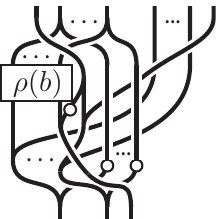}
\end{matrix}
\underset{\text{\bf Hf}}{\longleftrightarrow}
\\
\begin{matrix}
\includegraphics[scale=0.65]{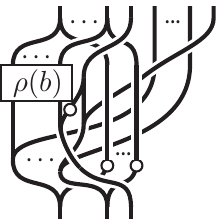}
\end{matrix}
=
\begin{matrix}
\includegraphics[scale=0.65]{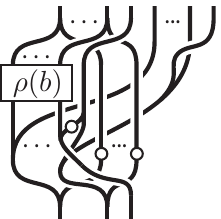}
\end{matrix}
=
\begin{matrix}
\includegraphics[scale=0.65]{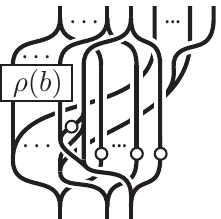}
\end{matrix}
\\
=
d(\sigma_n b)\circ(id^{\otimes (n+1)} \otimes f
\otimes 1)
\end{multline*}
\vspace{-8mm}
\caption{Two maps $d(\sigma_n b)$ and $d(\sigma_n b)\circ(id^{\otimes (n+1)} \otimes f\otimes 1)$ 
are equal as maps from $A^{\otimes 2(n+1)}$ to $A_{\sigma_nb}$.}
\label{fig:MII12}
\end{figure}
shows that 
$d(\sigma_n b)(\bx \otimes \by)-
d(\sigma_n b)(\bx \otimes f(\by)\otimes 1)\in I_{d(\sigma_n b)}$, 
and hence
$d(\sigma_n b)(\bx \otimes \by) \in I_{d(\sigma_n b)}$ since $\by \in \ker f$.  
Moreover, $d(\sigma_n b)(\bx \otimes \by) - (\varepsilon^{\otimes (n+1)})(\bx)\, \by \in I_{d(\sigma_n b)}$
by the definition of $I_{d(\sigma_n b)}$, 
we get
$(\varepsilon^{\otimes (n+1)})(\bx)\, \by \in I_{d(\sigma_n b)}$.  
Here $\bx$ is an arbitrary element of $A^{\otimes (n+1)}$, so $\by$ must be an element of $I_{d(\sigma_n b)}$.
\end{proof}
Next, we compare $A_b$ and $A_{\sigma_n^{-1}b}$.  
\begin{proposition}
For $b\in B_n$,  the $\Ad$-comodules $A_{b}$ and $A_{\sigma_n^{-1} b}$ are isomor\-phic.  
\label{prop:MII2}
\end{proposition}
\begin{proof}
Let $g$ be the linear surjection from $A^{\otimes (n+1)}$ to $A^{\otimes n}$ defined by 
\[
g = \mu_n\circ S_n^{-2}\circ \Psi_n^{-1}.
\]
\begin{figure}[htb!]
\begin{multline*}
g \circ d(\sigma_n^{-1} b)
=
\begin{matrix}
\includegraphics[scale=0.65]{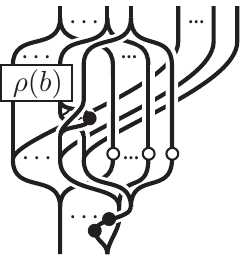}
\end{matrix}
=
\begin{matrix}
\includegraphics[scale=0.65]{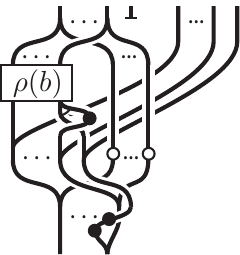}
\end{matrix}
=
\begin{matrix}
\includegraphics[scale=0.65]{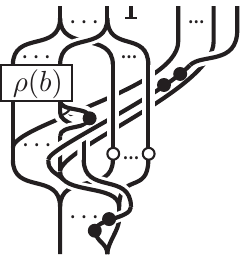}
\end{matrix}
\\
\underset{\text{(bc)}}{=}
\begin{matrix}
\includegraphics[scale=0.65]{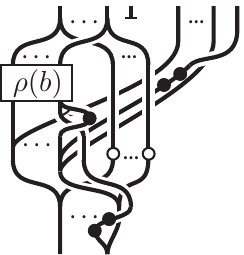}
\end{matrix}
=
\begin{matrix}
\includegraphics[scale=0.65]{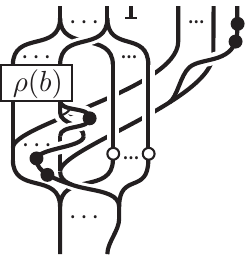}
\end{matrix}
=
\begin{matrix}
\includegraphics[scale=0.65]{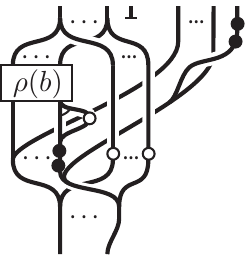}
\end{matrix}
\\
\underset{\text{Prop. \ref{prop:adS1}}}{=}
\begin{matrix}
\includegraphics[scale=0.65]{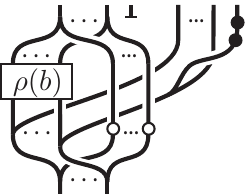}
\end{matrix}
=
d(b)\circ \varepsilon_n \circ(id^{\otimes (n+1)} \otimes g) 
\end{multline*}
\caption{The map $g\circ d(\sigma_n^{-1} b)$ is equal to $d(b)\circ \varepsilon_n \circ(id^{\otimes (n+1)} \otimes g)$. 
}
\label{fig:MII21}
\end{figure}
Then Figure \ref{fig:MII21} 
shows that $g\big(I_{d(\sigma_n^{-1} b)}\big) \subset I_{d(b)}$.  
To obtain the second equality in that Figure \ref{fig:MII21}, we slide down through
the $g$ part at the bottom to cancel the strand containing the right-most antipode, using the antipode axiom and anti-multiplicativity of $S^{-1}$.
The third equality similarly slides the top-rightmost strand through $g$.
\par
Next, Figure \ref{fig:MII22} 
shows that $\ker g \subset I_{d(\sigma_n^{-1} b)}$ using an argument similar to that of the 
previous proposition. 
 It follows that $g$ induces an isomorphism from 
$A_{\sigma_n^{-1} b}$ to $A_{b}$.  
\end{proof}
Instead of $I_{d(b)}$ we could also consider the $\Ad$-comodule $J_{d(b)}$ given by the image of
\begin{equation}
\boldsymbol{\mu}^{(n)}\circ\big(\boldsymbol{\Psi}^{(n)}\big)^{-1}\circ\Big(id^{\otimes n}\otimes \big(\rho(b) \circ\big(\boldsymbol{\theta}^{(n)}\big)^{-1}\circ \big(S^{-2}\big)^{\otimes n}\big)\Big) - \boldsymbol{\mu}^{(n)}
\label{eq:J}
\end{equation}
    where $\boldsymbol{\Psi}^{(n)}$ is the braiding of two bunches of $n$ strands and $\boldsymbol{\theta}^{(n)}$ is the full twist 
given in Figure \ref{fig:full}. The following proposition shows these are in fact isomorphic.
\begin{proposition}
$I_{d(b)}\cong J_{d(b)}$ as $\Ad$-comodules of $A^{\otimes n}$.
\end{proposition}
\begin{proof}
The deformation of Figure \ref{fig:proof1} shows that $J_{d(b)} \subset I_{d(b)}$.  
On the other hand, the deformation of Figure \ref{fig:proof2} shows that the image of 
$\boldsymbol{\mu}^{(n)} \circ 
(\boldsymbol{\mu}^{(n)} \otimes id^{\otimes n})\circ
\big(id^{\otimes n} \otimes \boldsymbol{\Psi}^{(n)}\big)\circ 
\big(b \otimes S^{\otimes n} \otimes id^{\otimes n}\big)\circ
\big(\boldsymbol{\Delta}^{(n)}\otimes id^{\otimes n}\big) - \varepsilon^{\otimes n} \otimes id^{\otimes n}$ is contained in $J_{d(b)}$.  
This means that $I_{d(b)} \subset J_{d(b)}$.  
\end{proof}
\begin{figure}[htb!]
\begin{multline*}
\begin{matrix}
\includegraphics[scale=0.6]{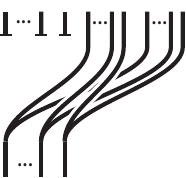}
\end{matrix}
\underset{\mod I_{d(b)}}{\equiv}
\begin{matrix}
\includegraphics[scale=0.6]{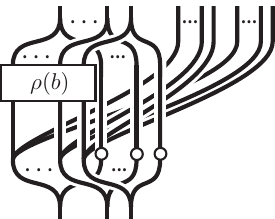}
\end{matrix}
=\ 
\begin{matrix}
\includegraphics[scale=0.6]{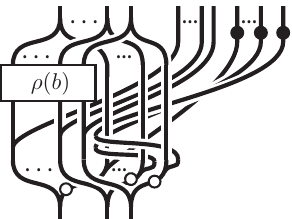}
\end{matrix}
\underset{{ \bf HmR' }}{=}
\\
\begin{matrix}
\includegraphics[scale=0.6]{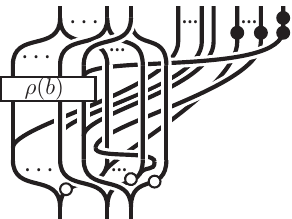}
\end{matrix}
\!\!\!\!
\underset{\text{Prop. \ref{prop:extramove}}}{=}\ 
\begin{matrix}
\includegraphics[scale=0.6]{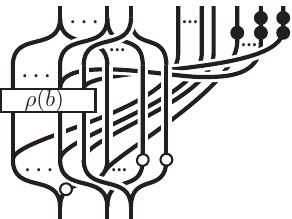}
\end{matrix}
\!\!\!\!
\underset{\text{Prop. \ref{prop:extramove}}}{=}\ 
\begin{matrix}
\includegraphics[scale=0.6]{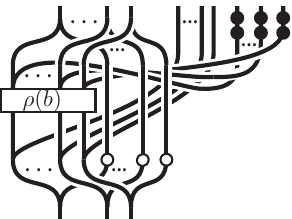}
\end{matrix}
\!\!\!\!
\underset{{ \text{Prop } \ref{prop:product} }}{=}
\\
\begin{matrix}
\includegraphics[scale=0.6]{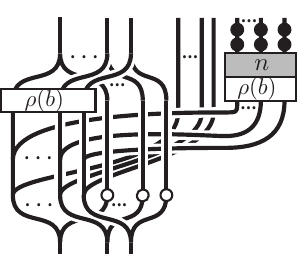}
\end{matrix}
\underset{\mod I_{d(b)}}{\equiv}
\begin{matrix}
\includegraphics[scale=0.6]{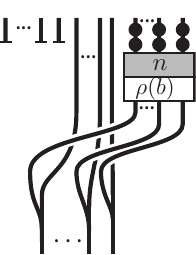}
\end{matrix},
\quad
\text{\footnotesize where $\begin{matrix}
\includegraphics[scale=0.6]{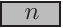}
\end{matrix} = \big(\boldsymbol{\theta}^{(n)}\big)^{-1}$.}
\end{multline*}
\vspace{-4mm}
\caption{Deformations to show that 
$\boldsymbol{\mu}^{(n)}(\bx\otimes \by)$  
is equal to 
$\boldsymbol{\mu}^{(n)}\circ\Big(\boldsymbol{\Psi}^{(n)}\Big)^{-1}\circ\Big(id^{\otimes n} \otimes \big(\rho(b) \circ \big(\boldsymbol{\theta}^{(n)}\big)^{-1} \circ S^{-2}\Big)(\bx\otimes \by)$ 
modulo 
$I_{d(b)}$.}
\label{fig:proof1}
\end{figure}
\begin{figure}[htb!]
\begin{multline*}
\begin{matrix}
\includegraphics[scale=0.6]{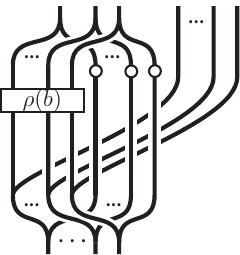}
\end{matrix}
\!\!\!\!
\underset{\mod J_b}{\equiv}
\begin{matrix}
\includegraphics[scale=0.6]{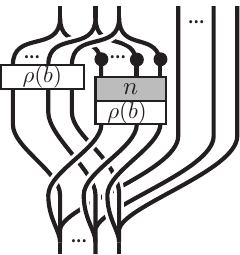}
\end{matrix}
\ =\ 
\begin{matrix}
\includegraphics[scale=0.6]{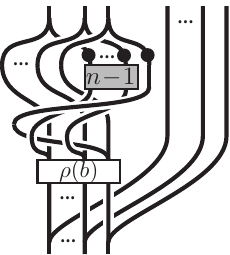}
\end{matrix}
\ =\ 
\begin{matrix}
\includegraphics[scale=0.6]{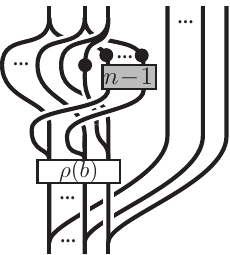}
\end{matrix}
\\
\ =
\begin{matrix}
\includegraphics[scale=0.6]{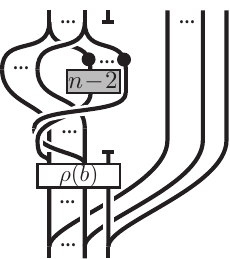}
\end{matrix}
\ =\ 
\cdots
\ =\ 
\begin{matrix}
\includegraphics[scale=0.6]{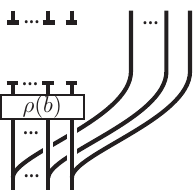}
\end{matrix}
\ =\ 
\begin{matrix}
\includegraphics[scale=0.6]{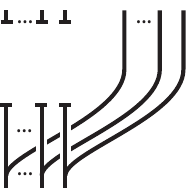}
\end{matrix}
\ =\ 
\begin{matrix}
\includegraphics[scale=0.6]{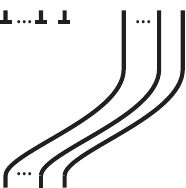}
\end{matrix}\ .
\end{multline*}
\vspace{-4mm}
\caption{Deformations to show that  $d(b)(\bx\otimes \by)$ is equal to $\big(\varepsilon^{\otimes n} \otimes id^{\otimes n}\big)(\bx\otimes \by)$ modulo $J_b$.
In the third equality we applied the antipode axiom to cancel the right-most antipode after shifting the box with the $n-1$ twist of the way to the far right.}
\label{fig:proof2}
\end{figure}
\subsection{Spanning set of $I_b$}
\begin{theorem}
Let $X$ be a set of generators of $A$ and
$x_i = 1^{\otimes(i-1)}\otimes x \otimes 1^{\otimes (n-i)}$
for $x \in X$.
Then the $\Ad$-comodule $I_{d(b)}$ in $A^{\otimes n}$ is spanned by 
\begin{equation}
\{d(b)(x_i \otimes \by)- {\varepsilon}(x)\by\,  \mid x \in X, \ i = 1, \cdots, n, \ \ \by \in A^{\otimes n}\}.  
\label{eq:idealgenerator}
\end{equation}
\label{thm:generator}
\end{theorem}
\begin{proof}
Let $I'$ be the $\Ad$-comodule spanned by 
$\{d(b)(x_i \otimes \by)- {\varepsilon}(x)\by\,  \mid x \in X, \ i = 1, \cdots, n, \ \ \by \in A^{\otimes n}\}$.
Since $I'\subset I_{d(b)}$ is obvious, we show that $I_{d(b)} \subset I'$.  
If $\big(d(b) - {\varepsilon}^{\otimes n}\otimes id^{\otimes n}\big)(\bx\otimes \by)$ and $\big(d(b) - {\varepsilon}^{\otimes n}\otimes id^{\otimes n}\big)(\bx' \otimes \by)$ are contained in $I'$ for any $\bx$, $\bx'$ and $\by$ in $A^{\otimes n}$,
then Figure \ref{fig:genproduct} shows that  $d(b) \big(\boldsymbol{\mu}^{(n)}(\bx \otimes \bx')\otimes \by\big)$ is equal to $\big(\varepsilon^{\otimes n}\big(\boldsymbol{\mu}^{(n)}(\bx \otimes \bx')\big) \by$ modulo $I'$.  
Hence $\big(d(b)-\varepsilon^{\otimes n}\otimes id^{\otimes n}\big) \big(\boldsymbol{\mu}^{(n)}(\bx \otimes \bx')\otimes \by\big)$ is contained in $I'$,
and this implies that  $I_{d(b)}\subset I'$.  
In the figure, $d(b)'$ means the part of $d(b)$ given by Figure \ref{fig:dbprime}.  
\begin{figure}[htb!]
\[
d({b}) =\ 
\begin{matrix}
\includegraphics[scale=0.6]{braidclosure0}
\end{matrix}
\ =
\ \ \begin{matrix}
\includegraphics[scale=0.6]{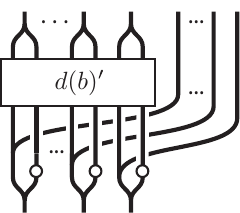}
\end{matrix}
\]
\caption{The part $d(b)'$ of the braided Hopf diagram $d(b)$.}
\label{fig:dbprime}
\end{figure}
\begin{figure}[htb!]
\begin{multline*}
d(b)\big(\boldsymbol{\mu}^{(n)}(\bx \otimes \bx')\otimes \by\big) =
\begin{matrix}
\includegraphics[scale=0.6]{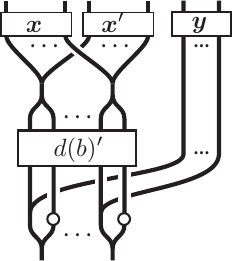}
\end{matrix}
\underset{\text{Prop. \ref{prop:product}}}{=}
\begin{matrix}
\includegraphics[scale=0.6]{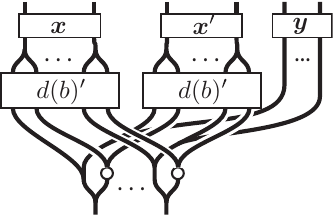}
\end{matrix}
\\
=
\begin{matrix}
\includegraphics[scale=0.6]{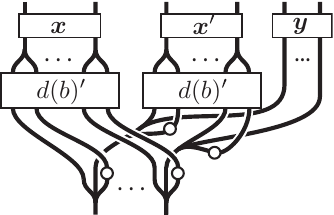}
\end{matrix}
\!\!\!\!
\underset{\mod I'}{\equiv}
\begin{matrix}
\includegraphics[scale=0.6]{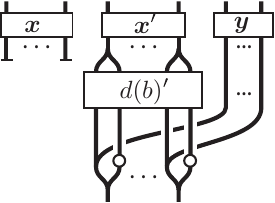}
\end{matrix}
\!\!\!\!
\underset{\mod I'}{\equiv}
\begin{matrix}
\includegraphics[scale=0.6]{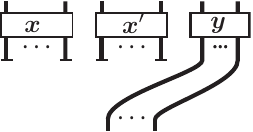}
\end{matrix}.
\end{multline*}
\caption{Computation of $d(b)\big(\boldsymbol{\mu}^{(n)}(\bx \otimes \bx')\otimes \by\big)$ modulo $I'$.}
\label{fig:genproduct} 
\end{figure}
\end{proof}
\section{Examples}\label{sec:examples}
Let $A$ be a finitely generated braided commutative braided Hopf algebra and let $X$ be a set of generators of $A$.  
We construct the space of  $A$ representations for the trivial knot, the Hopf link, the trefoil knot and the figure eight knot.  
\subsection{Trivial knot}
\begin{figure}[htb!]
\[
d(1) : 
\begin{matrix}
\includegraphics[scale=0.8]{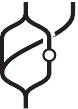}
\end{matrix}
=\ \ 
\begin{matrix}
\includegraphics[scale=0.8]{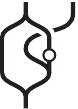}
\end{matrix}
=\ \ 
\begin{matrix}
\includegraphics[scale=0.8]{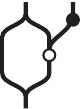}
\end{matrix}
\underset{\text{Prop. \ref{prop:extramove}}}{\longleftrightarrow}\ 
\begin{matrix}
\includegraphics[scale=0.8]{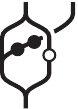}
\end{matrix}.
\]
\caption{Deformation of $d(1)$ to show $d(1)(x \otimes y)- d(1)\big(x \otimes S^{-2}(y)\big) \in I$.}
\label{fig:trivialdeform}
\end{figure}

Let $I=I_{d(1)}$ be the image of $d(1) - \varepsilon\otimes id$.
Then the space of $A$ representations for the trivial knot $A_1$ is given by $A_1 = A / I$.  
Since $d(1)$ can be deformed as Figure \ref{fig:trivialdeform}, we have 
$d(1)(x \otimes y)-  d(1)\big(x \otimes S^{-2}(y)\big)\in I$. 
We also have $d(1)(x \otimes y) - \varepsilon(x)y \in I$ and $d(1)\big(x \otimes S^{-2}(y)\big) - \varepsilon(x)S^{-2}(y) \in I$, hence $y - S^{-2}(y)\in I$.  
So, in the quotient space $A_1$, $S^2$ acts trivially.  
Moreover, from the relation \eqref{eq:J}, we have $\mu(x\otimes y) = \mu\circ \Psi(x\otimes  y)$.  
This implies that $A_1$ is a commutative algebra.  
\subsection{Hopf link}
\begin{figure}[htb!]
\begin{align*}
I_1 :& \quad  
\begin{matrix}
\includegraphics[scale=0.6]{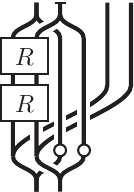}
\end{matrix}
=
\begin{matrix}
\includegraphics[scale=0.6]{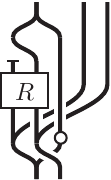}
\end{matrix}
=
\begin{matrix}
\includegraphics[scale=0.6]{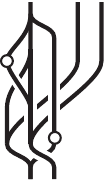}
\end{matrix},\\
I_2 :& \quad 
\begin{matrix}
\includegraphics[scale=0.6]{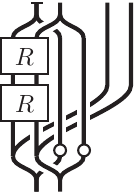}
\end{matrix}
\underset{\text{Prop. \ref{prop:product}}}{=}
\begin{matrix}
\includegraphics[scale=0.6]{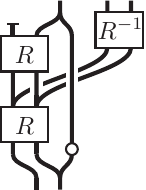}
\end{matrix}
\underset{\text{(bc)}}{=}
\begin{matrix}
\includegraphics[scale=0.6 ]{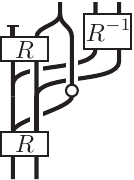}
\end{matrix}.
\end{align*}
\caption{The modules $I_1$ and $I_2$ for the Hopf link.}
\label{fig:hopfdeform}
\end{figure}

The Hopf link is the closure of $b=\sigma_1^2$ in $B_2$.  
Let $I_1 = \big(d(\sigma_1^2)-\varepsilon^{\otimes 2}\otimes id^{\otimes 2}\big)\big((A \otimes 1) \otimes A^{\otimes 2}\big)$ and 
$I_2 = \big(d(\sigma_1^2)-\varepsilon^{\otimes 2}\otimes id^{\otimes 2}\big)\big((1 \otimes A) \otimes A^{\otimes 2}\big)$.  
Then $I_{\sigma_1^2} = I_1 + I_2$ by Theorem \ref{thm:generator}.  
Figure \ref{fig:hopfdeform} shows that
$I_2 = \sigma_1 \cdot I_1$
and so we get 
$
A_{\sigma_1^2}
=
A \otimes A / (I_1 + \sigma_1\cdot I_1).  
$
\subsection{Trefoil knot}
\begin{figure}[htb!]
\[
I_1 : 
\begin{matrix}
\includegraphics[scale=0.6]{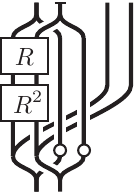}
\end{matrix}
=
\begin{matrix}
\includegraphics[scale=0.6]{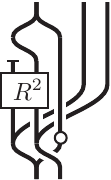}
\end{matrix},
\qquad
I_2 : 
\begin{matrix}
\includegraphics[scale=0.6]{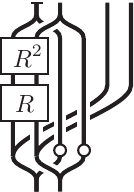}
\end{matrix}
\underset{\text{Prop. \ref{prop:product}}}{=}
\begin{matrix}
\includegraphics[scale=0.6]{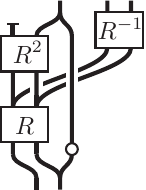}
\end{matrix}
\underset{\text{(bc)}}{=}
\begin{matrix}
\includegraphics[scale=0.6]{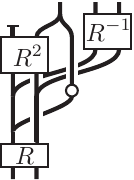}
\end{matrix}.
\]
\caption{The modules $I_1$ and $I_2$ for the trefoil knot.}
\label{fig:trefoildeform}
\end{figure}

The trefoil knot is the closure of $b=\sigma_1^3$ in $B_2$ so it can be treated like the Hopf link we considered above. 
In fact a similar computation
is valid for all closures of two strand braids.
Let $I_1 = \big(d(\sigma_1^3)-\varepsilon^{\otimes 2}\otimes id^{\otimes 2}\big)\big((A \otimes 1) \otimes A^{\otimes 2}\big)$ and 
$I_2 = \big(d(\sigma_1^3)-\varepsilon^{\otimes 2}\otimes id^{\otimes 2}\big)\big((1 \otimes A) \otimes A^{\otimes 2}\big)$.  
Then $I_{\sigma_1^2} = I_1 + I_2$ by Theorem \ref{thm:generator}.  
Figure \ref{fig:trefoildeform} shows that
$I_2 = \sigma_1 \cdot I_1$
and so we get 
$
A_{\sigma_1^3}
=
A \otimes A / (I_1 + \sigma_1\cdot I_1).  
$
\subsection{Figure eight knot}
The figure eight knot $4_1$ is isotopic to the closure of the braid $b =\sigma_2^{-1}\, \sigma_1\, \sigma_2^{-1}\, \sigma_1$.  
\begin{figure}[htb!]
\[
d(b) : 
\begin{matrix}
\includegraphics[scale=0.6]{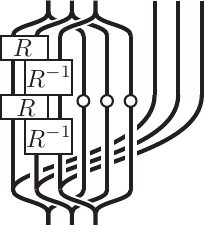}
\end{matrix}
\underset{\text{Prop. \ref{prop:product}}}{=}
\begin{matrix}
\includegraphics[scale=0.6]{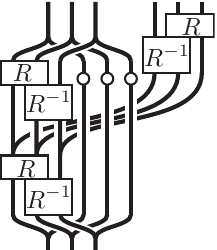}
\end{matrix}
\underset{\text{Prop. \ref{prop:product}}}{=}
\begin{matrix}
\includegraphics[scale=0.6]{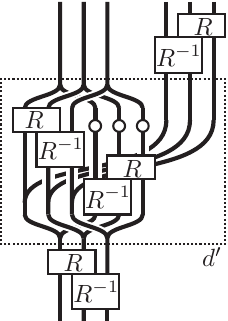}
\end{matrix}
\]
\caption{Deformation of $d(b)$ and the braided Hopf diagram $d'$.}
\label{fig:eightdeform0}
\end{figure}

Let $d'$ the braided Hopf diagram assigned in  Figure \ref{fig:eightdeform0} and
let $I'$ be the image of $d'  - \varepsilon^{\otimes 3} \otimes id^{\otimes 3}$.  
Then $\sigma_2^{-1}\sigma_1\,I'
=
\sigma_2^{-1}\sigma_1 \,\image(d'-  \varepsilon^{\otimes 3}\otimes  id^{\otimes 3})
=
\image\big((d(b)-  \varepsilon^{\otimes n} \otimes id^{\otimes 3})\circ(id^{\otimes 3} \otimes \sigma_2^{-1}\sigma_1)\big)
= I$
since $\sigma_2^{-1}\sigma_1$ is an automorphism of $A^{\otimes 3}$.  
Hence
$A_b$ is isomorphic to $A^{\otimes 3}/I'$.  
Let 
\[
I_1 = d'(A \otimes 1^{\otimes 2} \otimes A^{\otimes 3}), \ 
I_2 = d'(1 \otimes A\otimes 1 \otimes A^{\otimes 3}), \ 
I_3 = d'(1^{\otimes 2}\otimes A \otimes A^{\otimes 3}),
\]
then $I' = I_1 + I_2 + I_3$.  
\par

We first look at $I_3$.  
\begin{figure}[htb!]
\text{$I_3$ is spanned by $
\left(\begin{matrix}
\includegraphics[scale=0.6]{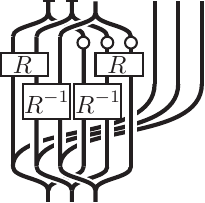}
\end{matrix}
-
\begin{matrix}
\includegraphics[scale=0.6]{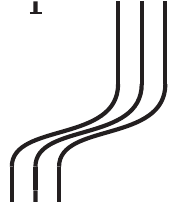}
\end{matrix}\right)$ where}
\begin{multline*}
\begin{matrix}
\includegraphics[scale=0.6]{eightdeform20}
\end{matrix}
=
\begin{matrix}
\includegraphics[scale=0.6]{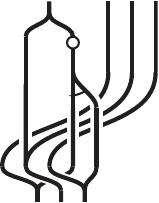}
\end{matrix}
\underset{\text{Prop. \ref{prop:extramove}}}{=}
\begin{matrix}
\includegraphics[scale=0.6]{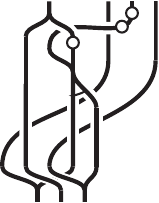}
\end{matrix}
=
\begin{matrix}
\includegraphics[scale=0.6]{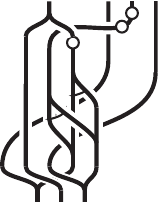}
\end{matrix}
=
\\
\begin{matrix}
\includegraphics[scale=0.6]{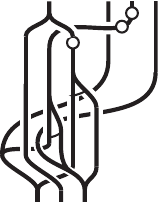}
\end{matrix}
=
\begin{matrix}
\includegraphics[scale=0.6]{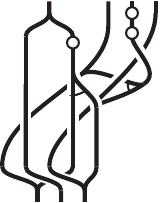}
\end{matrix}
\left(
\underset{\mod I'}{\equiv}
\begin{matrix}
\includegraphics[scale=0.6]{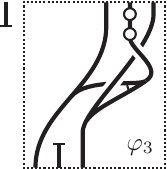}
\end{matrix}\ \ 
\right)
\end{multline*}
\caption{The subspace $I_3$.}
\label{fig:eightdeform3}
\end{figure}
Let $\varphi_3$ be the map from $A^{\otimes 3}$ to $A^{\otimes 2}$ defined by
\[
\varphi_3(x_1 \otimes x_2\otimes x_3) 
= 
\mu_1\circ \mu_3\circ \Psi_2 \circ 
(id^{\otimes 2}\otimes \ad) \circ \Psi_2\circ S_2^2
(x_1\otimes x_2 \otimes x_3).
\]
see Figure \ref{fig:eightdeform3}. The same figure shows that 
\begin{equation}
x_1 \otimes x_2 \otimes x_3- \varphi_3(x_1 \otimes x_2 \otimes x_3)\in I'.
\label{eq:I3}
\end{equation}  
Therefore, $A_b \cong A^{\otimes 2}/\big(\varphi_3(I_1) + \varphi_3(I_2) + \varphi_3(I_3)\big)$ where $A^{\otimes 2} = A \otimes \mathbb{C} \otimes A$.  
\par
\begin{figure}[htb!]
\[
\begin{matrix}
\\
\text{\footnotesize{$g_1 \ g_2 \ g_3$}}
\\
\includegraphics[scale=0.6]{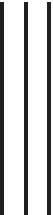}
\\
\text{\footnotesize{$g_1 \ g_2 \ g_3$}}
\end{matrix}
\underset{\mod I'}{\equiv}
\hspace{-1.5cm}
\begin{matrix}
\text{\footnotesize{$g_3^{-1}g_1^{}g_3^{}g_1^{-1}g_3^{}g_1^{}g_3^{-1}g_1^{-1}g_3^{}\qquad\qquad\qquad\qquad$}}
\\
\text{\footnotesize{$\qquad\quad g_3^{-1}g_1^{}g_3^{}$}}
\\
\text{\footnotesize{$\qquad \qquad g_3$}}
\\[-26pt]
\includegraphics[scale=0.6]{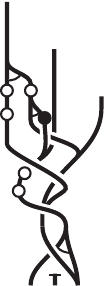}
\\
\text{\footnotesize{$g_1 \ g_2 \ g_3$}}
\end{matrix}
\hspace{-1cm}
\begin{tabular}{l}
Relations of $\pi_1(S^3 \setminus 4_1)$ obtained
\\
 from the picture:
\\{}\\
$g_1^{} = g_3^{-1}g_1^{}g_3^{}g_1^{-1}g_3^{}g_1^{}g_3^{-1}g_1^{-1}g_3^{}$
\\
$g_2 = g_3^{-1}g_1^{}g_3^{}$
\\
$g_3 = g_3$
\end{tabular}
\]
\caption{$x_1 \otimes  x_2 \otimes x_3 \equiv 
 \varphi_3 \big(1 \otimes \varphi_2(x_1 \otimes  x_2 \otimes x_3)\big) \mod I'$.}
 \label{fig:I23}
\end{figure}
\begin{figure}[htb!]
\text{$I_1$ is spanned by $
\left(\begin{matrix}
\includegraphics[scale=0.6]{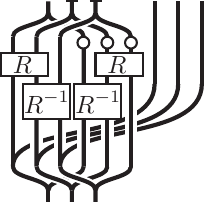}
\end{matrix}
-
\begin{matrix}
\includegraphics[scale=0.6]{eightdeform0}
\end{matrix}\right)
=
\left(
\begin{matrix}
\includegraphics[scale=0.6]{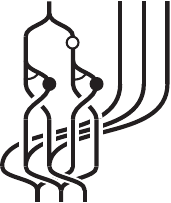}
\end{matrix}
-
\begin{matrix}
\includegraphics[scale=0.6]{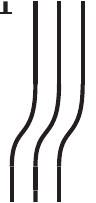}
\end{matrix}
\right)$.}
\caption{The subspace $I_1$.}
\label{fig:eightdeform1}
\end{figure}

Next, we look at $I_2$.  
\begin{figure}[htb!]
\text{$I_2$ is spanned by $
\left(\begin{matrix}
\includegraphics[scale=0.6]{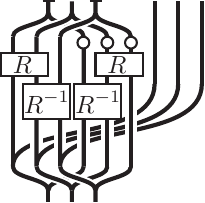}
\end{matrix}
-
\begin{matrix}
\includegraphics[scale=0.6]{eightdeform0}
\end{matrix}\right)$ where}
\begin{multline*}
\begin{matrix}
\includegraphics[scale=0.6]{eightdeform4}
\end{matrix}
=
\begin{matrix}
\includegraphics[scale=0.6]{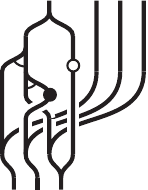}
\end{matrix}
\underset{\text{Prop. \ref{prop:admult}}}{=}
\begin{matrix}
\includegraphics[scale=0.6]{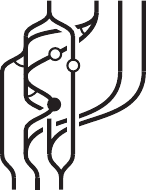}
\end{matrix}
\underset{\text{Prop. \ref{prop:extramove}}}{=}
\begin{matrix}
\includegraphics[scale=0.6]{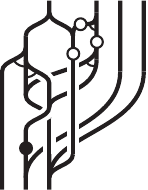}
\end{matrix}
\underset{\text{(bc)}}{=}
\\
\begin{matrix}
\includegraphics[scale=0.6]{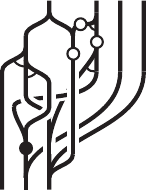}
\end{matrix}
=
\begin{matrix}
\includegraphics[scale=0.6]{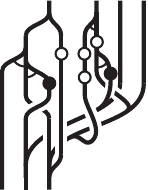}
\end{matrix}
\left(
\underset{\mod I'}{\equiv}
\begin{matrix}
\includegraphics[scale=0.6]{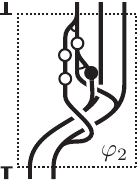}
\end{matrix}\ \ 
\right)
\end{multline*}
\caption{The subspace $I_2$. }
\label{fig:eightdeform2}
\end{figure}
Let  $\varphi_2$ be a map from $A^{\otimes 3}$ to $A^{\otimes 2}$ defined by
\begin{multline*}
\varphi_2(x_1 \otimes x_2\otimes x_3) 
= 
\mu_2 \circ \Psi_2 \circ \Psi_1 \circ
\mu_3\circ\mu_2 \circ \Psi_3\circ\Psi_2^{-1}\circ\mu_3 \circ
S_3^{-1}\circ
\\
 (id \otimes \ad\otimes id^{\otimes 2})\circ 
S_2 \circ S_1^2\circ(\ad\otimes id^{\otimes 2})
(x_1\otimes x_2 \otimes x_3)
\end{multline*}
as in Figure \ref{fig:eightdeform2}. 
The same figure shows that
$x_1 \otimes  x_2 \otimes x_3 - 1 \otimes \varphi_2(x_1 \otimes  x_2 \otimes x_3) \in I'$.  
Combining \eqref{eq:I3}, we get
\[
x_1 \otimes  x_2 \otimes x_3 - 
 \varphi_3 \big(1 \otimes \varphi_2(x_1 \otimes  x_2 \otimes x_3)\big) \in I'.  
\]
This relation is presented graphically in Figure \ref{fig:I23}.
Reading the diagram bottom to top and interpreting it in the group algebra we find the following presentation of 
$\pi_1(S^3\setminus 4_1)$.
\[
\pi(S^3\setminus 4_1) = 
\left<
g_1, g_3 \mid g_3^{-1}g_1^{}g_3^{}g_1^{-1}g_3^{}g_1^{}g_3^{-1}g_1^{-1}g_3^{}g_1^{-1}
\right>.
\]
%
Finally, the subspace $I_1$ is spanned by the image of the map given in Figure \ref{fig:eightdeform1}.  

\end{document}